\documentclass[11pt]{article}

\usepackage{amssymb}
\usepackage{amsfonts}
\usepackage{amsmath}
\usepackage{color}
\usepackage{multirow}

\newtheorem{theorem}{Theorem}

\newtheorem{corollary}[theorem]{Corollary}

\newtheorem{lemma}[theorem]{Lemma}

\renewcommand{\Re}{{\mathfrak{Re}}}

\begin{document}

\title{A two-type Bellman--Harris process initiated by a large number of
particles\thanks{ Supported by grant RFBR 11-01-00139 and by the
Programs of RAS ``Dynamical Systems and Control Theory''
(V.\,Vatutin) and ``Development of the methods of the
investigation of the stochastic models of those are oriented to
the population and biomedical applications''(V.\,Topchii)}}
\author{Vladimir Vatutin\thanks{
Steklov Mathematical institute RAS, Gubkin str. 8, Moscow, 119991,
Russia; e-mail: vatutin@mi.ras.ru}, \ \ Alexander Iksanov\thanks{
Faculty of Cybernetics, National Taras Shevchenko University of
Kyiv, 64/13 Volodymyrska, str., Kyiv-01601, Ukraine; e-mail:
iksan@univ.kiev.ua}, \ \  Valentin Topchii\thanks{ Sobolev
Institute of Mathematics SB RAS, Novosibirsk, 630090, Russia;
e-mail: topchij@ofim.oscsbras.ru}} \maketitle

\begin{abstract}
We investigate a two-type critical Bellman--Harris branching process with
the following properties: the tail of the life-length distribution of the
first type particles is of order $o(t^{-2})$; the tail of the life-length
distribution of the second type particles is regularly varying at infinity
with index $-\beta$, $\beta \in (0,1]$; at time $t=0$ the process starts
with a large number $N$ of the second type particles and no particles of the
first type. It is shown that the time axis $0\leq t<\infty $ splits into
several regions whose ranges depend on $\beta $ and the ratio $N/t$ within
each of which the process at time $t$ exhibits asymptotics (as $%
N,t\rightarrow \infty$) which is different from those in the other regions.
\end{abstract}

\section{Introduction}

To some extent, the present paper may be viewed as a continuation of \cite%
{VT13} in which the asymptotics of the survival probability of a two-type
critical Bellman--Harris branching process $\mathbf{Z}%
(t)=(Z_{1}(t),Z_{2}(t)) $, $t\geq 0$, was investigated and several
conditional limit theorems were proved for the distribution of the number of
particles at a distant time $t$ given that the process survives up to this
time. Since we will use some results obtained in the aforementioned paper we
recall some definitions and assertions given in \cite{VT13}.

The model of the two-type Bellman-Harris branching process in focus may be
informally described as follows. A particle of type $i\in \{1,2\}$ has the
life--length distribution $G_i(t)$, and at the end of her life she produces $%
\xi _{i1}$ particles of the first type and $\xi _{i2}$ particles of the
second type in accordance with generating function
\begin{equation*}
f_{i}(\mathbf{s})=f_{i}(s_{1},s_{2}):=\mathbf{E}\left[ s_{1}^{\xi
_{i1}}s_{2}^{\xi _{i2}}\right],\ \mathbf{s}:=(s_{1},s_{2})\in \left[ 0,1 %
\right] ^{2}.
\end{equation*}
Each particle ever born behaves in a similar manner, and she lives and
produces offspring independently of the co-existing particles and the past
history of the process.

Using the symbol $^{\dagger }$ for the transposition of vectors we introduce
two--dimensional vector-columns $\mathbf{G}(t):=(G_{1}(t),G_{2}(t))^{
\dagger }$, $\mathbf{f}(\mathbf{s})=\left( f_{1}(\mathbf{s}),f_{2}(\mathbf{s}
)\right) ^{\dagger }$. Symbols $\mathbf{1}$ and $\mathbf{0}$ will be used to
denote (depending on the context) either the vector-rows $(1,1)$ and $(0,0)$
or the vector-columns $(1,1)^{\dagger }$ and $(0,0)^{\dagger }$.

Put
\begin{equation}
m_{ij}:=\mathbb{E}\xi _{ij}=\left. \frac{\partial f_{i}(\mathbf{s})}{%
\partial s_{j}}\right\vert _{\mathbf{s}=\mathbf{1}},\ b_{jk}^{i}:=\mathbb{E}%
\xi _{ij}\xi _{ik}=\left. \frac{\partial ^{2}f_{i}(\mathbf{\ s})}{\partial
s_{j}\partial s_{k}}\right\vert _{\mathbf{s}=\mathbf{1}},\ i,j,k=1,2
\label{Fvarian}
\end{equation}%
and
\begin{equation}
\mathbf{M}(t):=\left(
\begin{array}{cc}
m_{11}G_{1}(t) & m_{12}G_{1}(t) \\
m_{21}G_{2}(t) & m_{22}G_{2}(t)%
\end{array}%
\right) ,\ \mathbf{M=M}(\infty ):=\left(
\begin{array}{cc}
m_{11} & m_{12} \\
m_{21} & m_{22}%
\end{array}%
\right) .  \label{01}
\end{equation}%
We assume that $\mathbf{Z}(t)$ is an indecomposable, aperiodic and
critical
process. This means, in particular, that there exists a positive integer $%
n_{0}$ such that all the elements of the matrix
$\mathbf{M}^{n_{0}}$ are  positive, the Perron root of
$\mathbf{M}$
is equal to 1 and there exist unique left and right eigenvectors $\mathbf{v}%
=(v_{1},v_{2})$ and $\mathbf{u}=(u_{1},u_{2})$ such that
\begin{equation}
\mathbf{Mu}^{\dagger }=\mathbf{u}^{\dagger },\,\mathbf{vM}=\mathbf{v},\,%
\mathbf{vu}^{\dagger }=1,\,\mathbf{u>0},\,\mathbf{v>0},\ \mathbf{v1}=1.
\label{DefMatrix}
\end{equation}%
In addition, we suppose that
\begin{equation}
B:=\frac{1}{2}\sum_{i,j,k=1,2}v_{i}b_{jk}^{i}u_{j}u_{k}<\infty .
\label{DefA}
\end{equation}

Observe that for the critical two-type indecomposable aperiodic process the
inequalities $m_{ii}<1$, $m_{11}+m_{22}>0$ and $m_{ij}>0$ for $i\neq j$ hold
true. Along with the criticality we impose the following conditions on the
tail behavior of the life-length distributions of particles:
\begin{equation}
1-G_{1}(t)=o(t^{-2})  \label{02}
\end{equation}
and
\begin{equation}
1-G_{2}(t)=\ell(t)t^{-\beta },\ \ \beta \in (0,1],  \label{03}
\end{equation}
where $\ell(t)$ is a function slowly varying at infinity. Here and hereafter
all unspecified limit relations are assumed to hold, as $t\to\infty$ or $%
N,t\to\infty$. Which of the alternatives prevails should be clear from the
context.

Let
\begin{equation}
\mu _{i}(t):=\int_{0}^{t}(1-G_{i}(w))\mathrm{d}w\ \text{and}\ \mu _{i}:=\mu
_{i}(\infty ).  \label{DefMu}
\end{equation}%
Clearly, $\mu _{1}<\infty $ and either $\mu _{2}<\infty $ and then $\beta =1$%
, or $\mu _{2}=\infty $ and then
\begin{equation}
\mu _{2}(t)\sim \left\{
\begin{array}{lll}
\frac{\textstyle{t^{1-\beta }\ell (t)}}{\textstyle{1-\beta }}=\frac{%
\textstyle t\left( 1-G_{2}(t)\right) }{\textstyle1-\beta }, & \text{if} &
\beta \in (0,1), \\
\ell _{1}(t), & \text{if} & \beta =1,%
\end{array}%
\right.   \label{AsymMu}
\end{equation}%
where $\ell _{1}(t):=\int_{0}^{t}\ell (u)u^{-1}\mathrm{d}u\rightarrow \infty
$ is a function slowly varying at infinity and $\ell (t)=o(\ell _{1}(t))$
(see, for example, Proposition 1.5.9a in \cite{BGT}). Observe that $\underset%
{t\rightarrow \infty }{\lim }\,\ell _{1}(t)=\mu _{2}$ irrespective of
whether $\mu _{2}$ is finite or not, and that $\mu _{i}(t)\leq t$. Thus,
\begin{equation}
1\leq R(t):=\frac{t}{\mu _{2}(t)}\sim \left\{
\begin{array}{lll}
\frac{\textstyle{(1-\beta )t^{\beta }}}{\textstyle{\ell (t)}}=\frac{%
\textstyle{1-\beta }}{\textstyle1-G_{2}(t)}, & \text{if } & \beta \in (0,1),
\\
\frac{\textstyle{t}}{\textstyle{\ell _{1}(t)}}, & \text{if} & \beta =1.%
\end{array}%
\right.   \label{012aa}
\end{equation}

Now we need more notation. For vectors $\mathbf{s}:=(s_{1},s_{2})\in \left[
0,1\right] ^{2}$ and $\mathbf{z}:=(z_{1},z_{2})\in \mathbb{Z}%
_{+}^{2}:=\{0,1,\ldots \}\times \{0,1,\ldots \}$ we write $\mathbf{s}^{%
\mathbf{z}}:=s_{1}^{z_{1}}s_{2}^{z_{2}}$ and for vectors $\mathbf{x}=\left(
x_{1},x_{2}\right) ^{\dagger }$ and $\mathbf{\ y}=\left( y_{1},y_{2}\right)
^{\dagger }$ we write
\begin{equation*}
\mathbf{x}\otimes \mathbf{y}:=\left( x_{1}y_{1},x_{2}y_{2}\right) ^{\dagger
}.
\end{equation*}%
Put further
\begin{equation*}
F_{i}(t;\mathbf{s})=F_{i}(t;s_{1},s_{2}):=\mathbf{E}_{i}\mathbf{s}^{\mathbf{%
\ Z}(t)},\quad \mathbf{F}(t;\mathbf{s}):=(F_{1}(t;\mathbf{s}),F_{2}(t;%
\mathbf{\ s}))^{\dagger },
\end{equation*}%
where here and hereafter
\begin{equation*}
\mathbf{E}_{j}\left[ \cdot \right] :=\mathbf{E}\left[ \cdot |\mathbf{Z}%
(0)=(\delta _{1j},\delta _{2j})\right] ,\ \mathbf{P}_{j}\left( \cdot \right)
:=\mathbf{P}\left( \cdot |\mathbf{Z}(0)=(\delta _{1j},\delta _{2j})\right)
,\ j=1,2,
\end{equation*}%
and $\delta _{ij}$ is the Kronecker delta.

The vector-function $\mathbf{F}(t;{\mathbf{s}})$ satisfies the following
system of integral equations (see, for instance, \cite{VT13})
\begin{equation}
\mathbf{F}(t;{\mathbf{s}})=\mathbf{s}\otimes (\mathbf{1}-\mathbf{G}
(t))+\int_0^t\mathbf{f}(\mathbf{F}(t-w;{\mathbf{s}}))\otimes \mathrm{d}
\mathbf{G }(w),  \label{f021}
\end{equation}
which, by setting $\mathbf{Q}(t;\mathbf{s}):=\mathbf{1}-\mathbf{F}(t,{\
\mathbf{s}})$, may be rewritten as
\begin{equation}
\mathbf{Q}(t;\mathbf{s})=\left( \mathbf{1}-\mathbf{s}\right) \otimes \left(
\mathbf{1}-\mathbf{G}(t)\right) +\int_{0}^{t}\left( \mathbf{1}-\mathbf{f}(
\mathbf{F}(t-w;\mathbf{s}))\right) \otimes \mathrm{d}\mathbf{G}(w).
\label{BBa1}
\end{equation}

Denote $\mathbf{\Phi }(\mathbf{s})=(\Phi _{1}(\mathbf{s}),\Phi _{2}(\mathbf{%
\ s }))^{\dagger }:=\mathbf{M}\mathbf{s}-(\mathbf{1}-\mathbf{f}(\mathbf{1}-
\mathbf{s}))$ and introduce a $2\times 2$ matrix
\begin{equation*}
\mathbf{G}_{\mathbf{I}}(t):=\left( G_{i}(t)\delta _{ij}\right) _{i,j=1}^{2}.
\end{equation*}
With this notation at hand a standard renewal argument enables us to rewrite %
\eqref{BBa1} in the form
\begin{eqnarray}
\mathbf{Q}(t;\mathbf{s}) &=&\left( \mathbf{1}-\mathbf{s}\right) \otimes
\left( \mathbf{1}-\mathbf{G}(t)\right) +\int_0^t \mathrm{d}\mathbf{M}(w)
\mathbf{Q} (t-w;\mathbf{s})  \notag \\
&&-\int_0^t \mathrm{d}\mathbf{G}_{\mathbf{I}}(w)\mathbf{\Phi }(\mathbf{Q}
(t-w; \mathbf{s})).  \label{EqvForQ}
\end{eqnarray}

One of the main characteristics of any critical branching process is its
survival probability. A specialization of Theorem 1 in \cite{13V} gives the
asymptotic behavior of the survival probability of the two-type
Bellman-Harris branching process which satisfies conditions \eqref{DefMatrix}
through \eqref{03}: for any fixed $\mathbf{s}:=(s_{1},s_{2})\in[0,1] ^{2}$
\begin{equation}
\mathbf{Q}(t;\mathbf{s})=\mathbf{1}-\mathbf{F}\left( t;\mathbf{s}\right)
\sim \mathbf{\ u}^{\dagger }\sqrt{\frac{v_{2}u_{2}}{B}\left(
1-G_{2}(t)\right) (1-s_{2})}.  \label{VatSurviv}
\end{equation}
In particular,
\begin{equation}
Q_{i}(t):=\mathbf{P}_{i}\left( \mathbf{Z}(t)\neq \mathbf{0}\right) \sim
\mathbf{P}_{i}\left( Z_{2}(t)>0\right) \sim u_{i}\sqrt{\frac{v_{2}u_{2}}{B}
\left( 1-G_{2}(t)\right) }  \label{12a}
\end{equation}
and, moreover,
\begin{equation}
\mathbf{P}_{i}\left( Z_{1}(t)>0\right) =o\left( Q_{i}(t)\right) .
\label{ASQ}
\end{equation}

The last two asymptotic relations mean that if the two-type population
survives up to a distant time $t$, then, with probability close to $1$, the
population at that time consists of the second type particles only. We
investigate this phenomenon in more detail in the situation when the
two-type Bellman-Harris process is initiated at time $t=0$ by a large number
$N$ of the second type particles and no particles of the first type, i.e., $%
\mathbf{Z}(0)=\left( Z_{1}(0),Z_{2}(0)\right) =\left( 0,N\right)$, and
analyze the distribution of the population size $\mathbf{Z} (t)=\left(
Z_{1}(t),Z_{2}(t)\right)$, as $t\rightarrow \infty$. Note that a similar
problem for a single-type critical Bellman-Harris branching process has been
investigated in \cite{V86}. The critical multitype Sevastyanov branching
processes (which are more general than the critical Bellman-Harris branching
processes, see \cite{Se}) initiated by a large number of particles were
(implicitly) considered in \cite{Shur76} under the assumption that the
expected life-lengths of particles of all types are finite. In view of (\ref%
{03}) our results do not follow (and, in fact, are essentially different)
from the results obtained in \cite{Shur76}.

Clearly,
\begin{equation}
\mathbf{E}\left[ \mathbf{s}^{\mathbf{Z}(t)}\Big|\mathbf{Z}(0)=\left(
0,N\right) \right] =F_{2}^{N}\left( t;\mathbf{s}\right) =e^{-N(1-F_{2}\left(
t;\mathbf{s }\right) )\left( 1+o(1)\right) }  \label{Blarge}
\end{equation}
provided that $\underset{t\to\infty}{\lim}\,(1-F_{2}\left( t;\mathbf{s}%
\right))=0$. Thus, to understand the asymptotic behavior of $\mathbf{Z}(t)$
under the present assumptions one has to investigate the behavior of $N(
1-F_{2}( t;\mathbf{s}))$, as $N,t\rightarrow \infty$, under a proper scaling
of the components of $\mathbf{Z}(t)$. If $N$ and $t$ tend to infinity in
such a way that $N\sqrt{1-G_{2}(t)}\rightarrow 0$, then, in view of (\ref%
{VatSurviv}) the population becomes extinct. If, however, $N\sqrt{
v_{2}u_{2}\left( 1-G_{2}(t)\right) /B}\rightarrow r\in (0,\infty )$, then
\begin{equation}
\lim_{N,t\rightarrow \infty }\mathbf{E}\left[ \mathbf{s}^{\mathbf{Z}(t)} %
\Big| \mathbf{Z}(0)=\left( 0,N\right) \right] =e^{-ru_{2}\sqrt{1-s_{2}}}.
\label{Z12}
\end{equation}
Thus, despite the indecomposability of the process there is only a
finite number of the second type particles and no particles of the
first type in the limit. This phenomenon has a natural intuitive
explanation: at a distant time $t$ the population only consists of
the particles whose life-length distributions have heavy tails
(see \cite{V86} where a similar effect is discussed for a
single-type critical Bellman-Harris process).

Below we list the basic assumptions of the paper.

\noindent \textbf{Hypothesis A.} \textit{The distribution function $G_{1}(t)$
satisfies \eqref{02} and the distribution function $G_{2}(t)$ satisfies %
\eqref{03}. In addition, if $\beta \in (0,1/2]$ then there exist positive
constants $C$ and $T_{0}$ such that for $t\geq T_{0}$ and any fixed $\Delta
>0$}
\begin{equation}
G_{2}(t+\Delta )-G_{2}(t)\leq C\Delta \ell(t)t^{-\beta -1}.  \label{UU2}
\end{equation}

Sometimes we will have to deal with ranges of $\beta $ other than $(0,1]$.
In these cases we write, say, that Hypothesis \textbf{A}$(a,b)$ or \textbf{A}%
$(a,b]$ holds, meaning that we \textbf{only} consider the range $\beta
\in(a,b)$ or $\beta \in (a,b]$ and require the validity of Hypothesis
\textbf{A} for the indicated range. In the sequel we denote by $C,\,C_{1},
\,C_{2},\,...$ positive constants whose specific values are of no
importance. The values of these constants need not be the same with each
usage.

Note that if Hypothesis \textbf{A} holds, then, according to Lemma \ref%
{L_firsrMom} below,
\begin{eqnarray}
N\mathbf{P}_{2}(Z_{1}(t) >0)=N(1-F_{2}(t;0,1)) \leq N\mathbf{E} _{2}\left[
Z_{1}(t)\right] \leq CN/\mu _{2}(t)  \label{OneNeglig}
\end{eqnarray}
which implies that there are no particles of the first type in the limit if $%
\mu _{2}(t)\gg N $. This means, in particular, that if, given $\mu_2(t)\gg N$%
, the limit
\begin{eqnarray}
\Pi_{2}(\lambda ) &:=&\lim_{N,t\rightarrow \infty }N\left( 1-\mathbf{E}_{2} %
\left[ e^{-\lambda Z_{2}(t)\psi (t)}\right] \right)  \notag \\
&=&\lim_{N,t\rightarrow \infty }N\left( 1-F_{2}\left(t;0,e^{-\lambda \psi
(t)}\right)\right)  \label{ome1}
\end{eqnarray}
exists for some function $\psi (t)$ and $\lambda >0$, then for any choice $%
s_1=s_1(t)\in \lbrack 0,1]$
\begin{equation}
\lim_{N,t\rightarrow \infty }N\left( 1-\mathbf{E}_{2}\left[
s_{1}^{Z_{1}(t)}e^{-\lambda Z_{2}(t)\psi (t)}\right] \right) =\Pi
_{2}(\lambda )  \label{ome2}
\end{equation}
and vice versa.

Here is a reasonable intuitive explanation of this effect. The branches of
the genealogical trees generated by $N$ initial particles consist of the
rays which may be thought of as those generated by renewal processes with
increments which (depending on the type of the corresponding particle) have
the distribution function $G_{1}(t)$ or $G_{2}(t)$. As we know by (\ref{Z12}%
), there are only a few surviving branches at a distant time $t$ such that $%
\mu _{2}(t)\gg N$ and, as a result, not too many rays attain the time-level $%
t$. Since the life-length distribution of the first type particles has a
light tail ($o(t^{-2})$), particles of this type are present in the
population at time $t$ with probability which is negligible in comparison
with the survival probability of the whole process up to this time.

It will be shown that there are several natural regions of $t=t(N)$ which
correspond to essentially different limiting distributions, as $%
N,t\rightarrow \infty $, of the vector $\mathbf{Z}(t),$ properly scaled.

The ranges of these regions depend on the behavior of the product $N(
1-G_{2}(t)),$ as $N,t\rightarrow \infty$, or, in view of the relation $%
\lim_{t\to\infty}R(t)( 1-G_{2}(t))=1-\beta$ which holds true whenever $\beta
\in \left( 0,1\right) $, on the behavior of the ratio $R(t)/N$. In the case $%
\beta=1$ the ranges of the regions will only be described in terms of $%
R(t)/N $. This motivates us to formulate the subsequent results in terms of $%
R(t)/N$ rather than in terms of a longer expression $N(1-G_{2}(t))$.

To describe the ranges of the regions in more detail we introduce three
functions $y= \mathfrak{g}_{1}(N)$, $y=\mathfrak{g}_{2}(N)$ and $y=\mathfrak{%
\ g}_{3}(N)$ which are the inverse functions to
\begin{equation*}
N(y)=(1-\beta )\frac{y^{\beta }}{\ell(y)},\ N(y)=\frac{y^{1-\beta }\ell(y) }{
1-\beta }\text{ and }N(y)=\sqrt{\frac{\textstyle{B}}{\textstyle{\
v_{2}u_{2}(1-G_{2}(y))}}},
\end{equation*}
respectively, if $\beta\in(0,1)$ and
\begin{equation*}
N(y)=\frac{y}{\ell_{1}(y)},\ N(y)=\ell_{1}(y)\text{ and }N(y)=\sqrt{\frac{ %
\textstyle{B}}{\textstyle{\ v_{2}u_{2}(1-G_{2}(y))}}},
\end{equation*}
respectively, if $\beta=1$. By Theorem 1.8.5 in \cite{BGT}
\begin{equation*}
\mathfrak{g}_{1}(N)=N^{1/\beta }L_{1}(N),\ \mathfrak{g}_{2}(N)=N^{1/\left(
1-\beta \right) }L_{2}(N),\ \mathfrak{g}_{3}(N)=N^{2/\beta }L_{3}(N)
\end{equation*}
for some functions $L_{i}(\cdot),$ $i=1,2,3$, slowly varying at infinity,
where $\beta\in (0,1]$, excluding the case $\beta=1$ for $\mathfrak{g}_{2}$.
Of course, $\mathfrak{g}_{2}(N)\gg N^{k}$ for all $k\in \mathbb{N}$, if $%
\beta=1$.

It can be checked that
\begin{equation*}
\mathfrak{g}_{2}(N)\ll \mathfrak{g}_{1}(N)\ll \mathfrak{g}_{3}(N)\text{ \
for \ }\beta \in \lbrack 0,1/2),
\end{equation*}
\begin{equation*}
\ \ \ \mathfrak{g}_{1}(N)\ll \mathfrak{g}_{2}(N)\ll \mathfrak{g}_{3}(N)\text{
\ for \ }\beta \in (1/2,2/3),
\end{equation*}
and
\begin{equation*}
\mathfrak{g}_{1}(N)\ll \mathfrak{g}_{3}(N)\ll \mathfrak{g}_{2}(N)\text{ \
for \ }\beta \in (2/3,1],
\end{equation*}
as $N\to\infty$. We write $t=t(N)\in \left( \mathfrak{g} _{i}(N),\mathfrak{g}
_{j}(N)\right)$ (or $\in \left( \mathfrak{g}_{i}(N), \mathfrak{g}
_{j}(N)\right) $) if
\begin{equation*}
\mathfrak{g}_{i}(N)\ll t(N)\ll \mathfrak{g}_{j}(N)
\end{equation*}
and $t=t(N)\simeq \mathfrak{g}_{i}(N)$ (or $\simeq \mathfrak{g}_{i}(N)$) if
\begin{equation*}
\lim_{N\rightarrow \infty }\frac{t(N)}{\mathfrak{g}_{i}(N)}=c\in \left(
0,\infty \right) .
\end{equation*}

It follows from the definitions above and properties of the functions $R(t)$
and $\mu _{2}(t)$ that
\begin{eqnarray*}
R(t)/N \rightarrow 0\ \, &\Leftrightarrow& t=o(\mathfrak{g}_{1}(N)), \\
R(t)/N \rightarrow r\ \, &\Leftrightarrow& t\sim \mathfrak{g}_{1}(Nr), \\
R(t)/N \rightarrow \infty &\Leftrightarrow& t\gg \mathfrak{g}_{1}(N)
\end{eqnarray*}
and
\begin{eqnarray*}
\mu _{2}(t)/N \rightarrow 0\ \, &\Leftrightarrow& t=o(\mathfrak{g}_{2}(N)),
\\
\mu _{2}(t)/N \rightarrow r_{\ast }\,&\Leftrightarrow& t\sim \mathfrak{g}
_{2}(Nr_{\ast }), \\
\mu _{2}(t)/N \rightarrow \infty &\Leftrightarrow& t\gg \mathfrak{g}_{2}(N).
\end{eqnarray*}

In the case $\beta\in (0,1)$ we will call the ranges of $t=t(N)$ satisfying,
as $N\rightarrow\infty$, either of the conditions
\begin{equation*}
N\left( 1-G_{2}(t)\right) \rightarrow \infty \ \ \text{ or } \ \
R(t)/N\rightarrow 0 \ \ \text{ or } \ \ t=o(\mathfrak{g}_{1}(N))
\end{equation*}
the \textit{early evolutionary stages} of the population, the ranges
satisfying either of the conditions
\begin{equation*}
N\left( 1-G_{2}(t)\right) \rightarrow r_{1}\in \left( 0,\infty \right) \ \
\text{ or } \ \ R(t)/N\rightarrow r\in \left( 0,\infty \right) \ \ \text{ or
} \ \ t\sim \mathfrak{g}_{1}(Nr)
\end{equation*}
the \textit{intermediate evolutionary stages}, and the ranges satisfying
either of the conditions
\begin{equation*}
N\left( 1-G_{2}(t)\right) \rightarrow 0 \ \ \text{ or } \ \
R(t)/N\rightarrow \infty \ \ \text{ or } \ \ t\gg \mathfrak{g}_{1}(N)
\end{equation*}
the \textit{final evolutionary stages}. In the case $\beta=1$ we will use
the same definition with the first condition (that involves $N\left(
1-G_{2}(t)\right)$) omitted.

Two remarks are in order. 1) Although the assumptions arising in the present
paper could have been formulated in terms of the time parameter $t=t(N)$ we
prefer to state them in terms of $R(t)/N$ and $\mu_2(t)/N$.

\noindent 2) The ensuing presentation will make it clear that the behavior
of the ratio $\mu _{2}(t)/N $, as $N,t\rightarrow \infty $, governs further
splitting the evolutionary stages just introduced into subranges which are
characterized by different limiting distributions of the vector $\mathbf{Z}%
(t)$, properly scaled.

The remaining part of the paper is structured as follows. While
the main results are stated in Section \ref{Sec1}, their proofs
are given in Section \ref{Sec3} (the early evolutionary stages),
Section \ref{Sec4} (the intermediate evolutionary stages) and
Section \ref{Sec5} (the final evolutionary stages). In Section
\ref{Sec2} we recall some known results taken mainly from
\cite{VT13}. These concern various asymptotic properties of
renewal matrices and generating functions, and are extensively
used throughout the paper.

\section{Main results \label{Sec1}}

\subsection{Early evolutionary stages}

In this subsection we investigate the asymptotic behavior of the number of
particles for the early evolutionary stages.

Denote by $\mathbf{D}=\left( D_{ij}\right) _{i,j=1}^{2}$ a $2\times 2$
matrix with positive entries
\begin{equation}
D_{ii}:=(1-m_{ii})\mu_2D \ \text{and} \ D_{ij}:=m_{ij}\mu_2D ,\quad i\neq j ,
\label{Dfin1}
\end{equation}
in the case when $\mu_2<\infty$, where $D:=((1-m_{22})\mu
_{1}+(1-m_{11})\mu_2)^{-1}$, and
\begin{equation}
D_{ii}=\frac{1-m_{ii}}{1-m_{11}}\text{ and }D_{ij}=\frac{m_{ij}}{1-m_{11}}
,\quad i\neq j,  \label{D1}
\end{equation}
in the case when $\mu_2=\infty$.

Set $\Gamma _\beta = 1$ for $\beta=1$ and
\begin{equation*}
\Gamma_\beta:=\frac{\sin {\pi \beta }}{\pi \beta (1-\beta )}
\end{equation*}
for $\beta \in (0,1)$.

\begin{theorem}
\label{T_initial1} Suppose that Hypothesis \textbf{A} holds and that
\begin{equation}  \label{inter13}
\lim_{N,\,t\rightarrow \infty }R(t)N^{-1}=0 \ \ \text{and} \ \
\lim_{N,\,t\rightarrow \infty }\mu_2(t)N^{-1}=0.
\end{equation}
If $\beta =1/2$, assume additionally that
\begin{equation}
\lim_{N,\,t\rightarrow \infty }\frac{\mu_{2}(t)}{N}\int_0^t\frac{\mathrm{d}w
}{(1+\mu _2(w))^{2}}=0.  \label{Redu}
\end{equation}
Then, for any $\lambda_1,\lambda_2>0$,
\begin{eqnarray*}
&&\lim_{N,\,t\rightarrow \infty }\mathbf{E}\left[ \exp \left\{ -\lambda _{1}
\frac{Z_{1}(t)\mu _{2}(t)}{N}-\lambda _{2}\frac{Z_{2}(t)}{N}\right\} \bigg|
\mathbf{Z}(0)=\left( 0,N\right) \right] \\
&&\qquad \qquad \qquad \qquad \qquad\quad =\exp \left\{ -\mu
_{1}\beta \Gamma _{\beta }D_{21}\lambda _{1}-D_{22}\lambda
_{2}\right\} .
\end{eqnarray*}
\end{theorem}

\begin{corollary}
\label{C_initial1} Suppose that Hypothesis \textbf{A} holds and that $%
\lim\limits_{N,\,t\rightarrow \infty }R(t)N^{-1}~=~0$. Then, for $\lambda >0$,
\begin{equation*}
\lim_{N,t\rightarrow \infty }\mathbf{E}\left[ \exp \left\{ -\lambda \frac{
Z_{2}(t)}{N}\right\} \bigg|\mathbf{Z}(0)=\left( 0,N\right) \right] =\exp
\left\{-D_{22}\lambda \right\} .
\end{equation*}
\end{corollary}

Note, that this corollary does not require $\lim_{N,\,t\rightarrow \infty
}\mu_2(t)N^{-1}=0$, nor condition (\ref{Redu}).

Set
\begin{equation*}
\mathbf{O}(s):=\left( O_{1}(s),\ O_{2}(s)\right) ^{\dag }:=\beta \Gamma
_{\beta }\int_{0}^{\infty }\mathbf{D\Phi }(\mathbf{Q}(w;s,1))\mathrm{d}w.
\end{equation*}

\begin{theorem}
\label{T_initial2} Suppose that Hypothesis \textbf{A}$(0,0.5]$ holds and
that
\begin{equation*}
\lim_{N,\,t\rightarrow \infty }R(t)N^{-1}=0 \ \ \text{and} \ \
\lim_{N,\,t\rightarrow \infty }\mu_2(t)N^{-1}= r^{-1}\in \left( 0,\infty
\right).
\end{equation*}
If $\beta =1/2$, assume additionally that
\begin{equation}
\Upsilon :=\int_0^\infty \frac{\mathrm{d}w}{(1+\mu _2(w))^{2}}<\infty.
\label{Sup_cond}
\end{equation}
Then, for any $s\in \lbrack 0,1]$ and $\lambda>0,$
\begin{eqnarray*}
&&\lim_{N,t\rightarrow \infty }\mathbf{E}\left[ s^{Z_{1}(t)}\exp \left\{
-\lambda \frac{Z_{2}(t)}{N}\right\} \bigg|\mathbf{Z}(0)=\left( 0,N\right) %
\right] \\
&&\qquad\qquad\qquad=\exp \left\{ -r\beta \Gamma _{\beta }\mu
_{1}D_{21}\left( 1-s\right) +rO_{2}(s)-D_{22}\lambda \right\} .
\end{eqnarray*}
\end{theorem}

Thus, in this case we asymptotically have a few individuals of the first
type, while the number of individuals of the second type is still of order~$%
N $. Moreover, $Z_1(t)$ and $Z_2(t)N^{-1}$ are asymptotically independent.

\noindent \textbf{Remark.} In view of (\ref{AsymMu}) and (\ref{012aa}) the
assumptions of Theorem \ref{T_initial2} hold if either $\beta =1/2$ and $%
\lim_{t\to\infty}\ell(t)=\infty$ or $\beta<1/2$.

\begin{theorem}
\label{T_initial3} Suppose that Hypothesis \textbf{A}$(0,0.5]$ holds and
that
\begin{equation*}
\lim_{N,t\rightarrow \infty}R(t)N^{-1}=0 \ \ \text{and} \ \
\lim_{N,t\rightarrow \infty }\mu _2(t)N^{-1}=\infty.
\end{equation*}
Then, for any $\lambda >0,$
\begin{eqnarray}
&&\hspace{-0.7cm} \lim_{N,t\rightarrow \infty }\mathbf{E}\left[ \exp \left\{
-\lambda \frac{Z_{2}(t)}{N}\right\} ;Z_{1}(t)=0\bigg|\mathbf{Z}(0)=\left(
0,N\right) \right]  \notag \\
&&\ =\lim_{N,t\rightarrow \infty }\mathbf{E}\left[ \exp \left\{ -\lambda
\frac{Z_{2}(t)}{N}\right\} \bigg|\mathbf{Z}(0)=\left( 0,N\right) \right]
=\exp \left\{-D_{22}\lambda \right\}.  \label{inter14}
\end{eqnarray}
\end{theorem}

\noindent \textbf{Remark.} In view of (\ref{AsymMu}) and (\ref{012aa}) the
assumptions of Theorem \ref{T_initial3} hold for $\beta =1/2$ only if $%
\lim_{t\to\infty}\ell(t)=\infty$.

\subsection{The intermediate evolutionary stages}

In this section we formulate theorems which describe the limiting behavior
of the population for the intermediate evolutionary ranges, i.e., we assume
that the limit of $R(t)/N$ is positive and finite. Unlike the early
evolutionary stages the asymptotic results here are affected by \textit{\
genuine properties} of branching processes. There are three essentially
different intermediate subranges which are characterized by one of the
conditions $\mu _{2}(t)N^{-1}\rightarrow \infty $, $\mu
_{2}(t)N^{-1}\rightarrow 0$ or $\mu _{2}(t)N^{-1}\rightarrow r_{1}\in \left(
0,\infty \right)$ which is assumed to hold along with the defining property
of the intermediate stages. We only analyze the first and the second
subranges. The remaining case $R(t)N^{-1}\rightarrow r\in \left( 0,\infty
\right)$ and $\mu _{2}(t)N^{-1}\rightarrow r_{2}\in \left( 0,\infty \right)$
which implies $\beta=1/2$ will be considered in a separate paper, for it
requires much more delicate analysis.

Put
\begin{equation}
N_{i}(\mathbf{x}):=\frac{1}{2}\sum_{j,k=1}^{2}b_{jk}^{i}x_{j}x_{k},\quad
\mathbf{N}\left( \mathbf{x}\right) :=\left( N_{1}(\mathbf{x}),N_{2}(\mathbf{%
\ x })\right) ^{\dagger }  \label{DefN}
\end{equation}
and observe that by Taylor's formula
\begin{eqnarray*}
&&\sum_{k=1}^{2}m_{ik}\left( 1-s_{k}\right) -\left( 1-f_{i}\left( \mathbf{s}
\right) \right) =\frac{1}{2}\sum_{j,k=1}^{2}b_{jk}^{i}\left( 1-s_{j}\right)
\left( 1-s_{k}\right) \\
&&+o\left( \left( 1-s_{1}\right) ^{2}+\left( 1-s_{2}\right) ^{2}\right)
=N_{i}(\mathbf{1}-\mathbf{s})+o\left( \left\Vert \mathbf{1}-\mathbf{s}
\right\Vert ^{2}\right)
\end{eqnarray*}
as $s_{1},s_{2}\uparrow 1$, or, equivalently,
\begin{equation}
\mathbf{\Phi }(\mathbf{1}-\mathbf{s})=\mathbf{M}\left( \mathbf{1-\mathbf{s}}
\right) \mathbf{-}\left( \mathbf{1}-\mathbf{f}(\mathbf{s})\right) =\mathbf{N}
(\mathbf{1}-\mathbf{s})+o\left( \mathbf{1}\left\Vert \mathbf{1}-\mathbf{s}
\right\Vert ^{2}\right) .  \label{DDEF}
\end{equation}
Letting
\begin{equation*}
\bar{b}:=\max_{i,j,k}b_{jk}^{i}
\end{equation*}
we infer
\begin{equation}
\mathbf{N}(\mathbf{x})\leq \bar{b}\left\Vert \mathbf{x}\right\Vert ^{2}
\mathbf{1}  \label{OcN}
\end{equation}
and
\begin{equation}
\left\Vert \mathbf{N}(\mathbf{x})-\mathbf{N}(\mathbf{x}^{\ast })\right\Vert
\leq \bar{b}\left( \left\Vert \mathbf{x}\right\Vert +\left\Vert \mathbf{x}
^{\ast }\right\Vert \right) \left\Vert \mathbf{x}-\mathbf{x}^{\ast
}\right\Vert .  \label{NorN}
\end{equation}

In the proof of Theorem \ref{T_Zone3} it will be shown that the system of
equations
\begin{equation}
\mathbf{\Omega }\left( \lambda \right) =\mathbf{D}\left( 0,\lambda \right)
^{\dagger }-\Gamma _{\beta }\int_{0}^{1}\frac{\mathbf{DN}\left( \mathbf{\
\Omega }\left( \lambda \left( 1-w\right) ^{\beta }\right) \right) }{\left(
1-w\right) ^{2\beta }}\mathrm{d}w^{\beta }, \ \ \lambda>0,  \label{UniNew1}
\end{equation}
has a unique solution with non-negative components, and we denote this
solution by $\mathbf{\Omega }\left( \lambda \right) :=\left( \Omega
_{1}(\lambda ),\Omega _{2}(\lambda )\right) ^{\dagger }$.

\begin{theorem}
\label{T_Zone3} Suppose that Hypothesis \textbf{A} holds and that
\begin{equation}
\lim_{N,\,t\rightarrow \infty }R(t)N^{-1}=1/r\in \left( 0,\infty \right).
\label{Rcons}
\end{equation}
Then, for any $\lambda >0,$
\begin{equation}
\lim_{N,\,t\rightarrow \infty }\mathbf{E}\left[ \exp \left\{ -\lambda \frac{
rZ_2(t)}{N}\right\} \bigg|\mathbf{Z}(0)=\left( 0,N\right) \right]
=e^{-r\Omega _{2}(\lambda)}.  \label{RconsGlobal}
\end{equation}

Furthermore, if Hypothesis \textbf{A}$(0,0.5]$ and condition $(\ref{Rcons})$
hold, and \newline
$\lim_{N,\,t\rightarrow \infty }\mu _{2}(t)N^{-1}=\infty $, then, for any $%
\lambda >0$,
\begin{equation}
\lim_{N,\,t\rightarrow \infty }\mathbf{E}\left[ \exp \left\{ -\lambda \frac{
rZ_{2}(t)}{N}\right\} ;Z_{1}(t)=0\bigg|\mathbf{Z}(0)=(0,N)\right]
=e^{-r\Omega _{2}(\lambda )}.  \label{Cr}
\end{equation}
\end{theorem}

Observe that there are no first type particles in the limit under the
asymptotic regime $\mu _{2}(t)N^{-1}\rightarrow \infty$. Note also that the
assumptions $R(t)N^{-1}\rightarrow r^{-1}$ and $\mu _{2}(t)N^{-1}\rightarrow
\infty $ entail $\beta =1/2$ and $\lim_{t\to\infty}\ell(t)=\infty$, or $%
\beta <1/2$.

Now we consider the case $R(t)N^{-1}\rightarrow r^{-1}$ and $%
\mu_{2}(t)N^{-1}\rightarrow 0$ which entails $\beta =1/2$ and $%
\lim_{t\to\infty} \ell(t)=0$, or $\beta>1/2$.

Put
\begin{equation*}
\mathbf{C}_{\beta }:=\left(
\begin{array}{cc}
D_{11}\beta \Gamma _{\beta } & D_{12} \\
D_{21}\beta \Gamma _{\beta } & D_{22}%
\end{array}
\right).
\end{equation*}
In the proof of Theorem \ref{T_Betamore} it will be shown that the system of
equations
\begin{equation}
\mathbf{H}(\theta ,\lambda )=\mathbf{C}_{\beta }\left( 1,\lambda \theta
^{1-\beta }\right) ^{\dagger }-\Gamma _{\beta }\theta ^{2\beta
-1}\int_{0}^{1}\frac{\mathbf{DN}\left( \mathbf{H}(\theta (1-y),\lambda
)\right) }{\left( 1-y\right) ^{2-2\beta }}\mathrm{d}y^{\beta },
\label{Systnew}
\end{equation}
for $\theta,\lambda>0$, with $\beta\in (1/2, 1]$, has a unique solution with
non-negative components, and we denote this solution by $\mathbf{H}(\theta
,\lambda )=\left( H_{1}(\theta ,\lambda ),H_{2}(\theta ,\lambda )\right)
^{\dagger }$.

\begin{theorem}
\label{T_Betamore} Suppose that Hypothesis \textbf{A}$(0.5,1]$ holds and
that
\begin{equation*}
\lim_{N,\,t\to\infty} R(t)N^{-1}=r^{-1}\in (0,\infty) \ \ \text{and} \ \
\lim_{N,\,t\to\infty} \mu _{2}(t)N^{-1}=0.
\end{equation*}
Then, for $\lambda _{1}>0$ and $\lambda _{2}>0$,
\begin{eqnarray*}
&&\lim_{N,t\rightarrow \infty }\mathbf{E}\left[ \exp \left\{ -\lambda _{1}
\frac{\mu _{2}(t)Z_{1}(t)}{\mu _{1}R(t)}-\lambda _{2}\frac{Z_{2}(t)}{R(t)}
\right\} \bigg|\mathbf{Z}(0)=(0,N)\right] \\
&&\qquad \qquad \qquad \qquad \qquad =\exp \left\{ -r\lambda _{1}H_{2}\left(
\lambda _{1}^{1/\left( 2\beta -1\right) },\lambda _{2}\lambda _{1}^{-\beta
/\left( 2\beta -1\right) }\right) \right\} .
\end{eqnarray*}
\end{theorem}

\subsection{Final evolutionary stages}

In the first part of Theorem \ref{T_zone2} given next we investigate the
asymptotic behavior of the number of particles under the conditions
\begin{equation}
R(t)N^{-1}\rightarrow \infty \ \text{and} \ N\sqrt{1-G_{2}(t)}\rightarrow
\infty,  \label{FinFirst}
\end{equation}
the first of these being a defining property of the final evolutionary
stages. In the second part of that theorem we work under the conditions %
\eqref{FinFirst} and
\begin{equation}
N/\mu _{2}(t)\rightarrow 0.  \label{FinFirst1}
\end{equation}
Relations \eqref{FinFirst} and \eqref{FinFirst1} together imply
\begin{equation}
0\longleftarrow \frac{N}{\mu _{2}(t)}=\frac{N\left( 1-G_{2}(t)\right) ^{1/2}
}{\mu _{2}(t)\left( 1-G_{2}(t)\right) ^{1/2}}=O\bigg( \frac{N\left(
1-G_{2}(t)\right) ^{1/2}}{t\left( 1-G_{2}(t)\right) ^{3/2}}\bigg)
\label{RatNm}
\end{equation}
and thereupon
\begin{equation}
\lim_{t\to\infty} t\left( 1-G_{2}(t)\right) ^{3/2}=\infty  \label{LimG23}
\end{equation}
which means that conditions (\ref{FinFirst}) and \eqref{FinFirst1} may only
hold simultaneously if either $\beta =2/3$ and $\lim_{t\to\infty}\ell(t)=
\infty$ or $\beta <2/3$.

\begin{theorem}
\label{T_zone2} Suppose that Hypothesis \textbf{A}$(0,1)$ holds, and that
\begin{equation}
\lim\limits_{N,\,t\rightarrow \infty }\frac{N\psi (t)}{ \mu_2(t)}=0,
\label{M1}
\end{equation}
\begin{equation}
1/R(t)=o(\psi (t)),  \label{XX}
\end{equation}
and
\begin{equation}
\lim_{N,\,t\to\infty} N\sqrt{\frac{v_{2}u_{2}}{B} \psi (t)\left(
1-G_{2}(t)\right)}=r\in\left(0,\infty\right)  \label{Nzone1}
\end{equation}
for a non-increasing, regularly varying function $\psi (t):= t^{-\gamma
}\ell_2(t)$, $\gamma \in [0,1)$, such that $\lim_{t\to\infty}\psi(t)=0$.
Then, for each $\lambda>0$,
\begin{equation}
\lim_{N,t\rightarrow \infty }\mathbf{E}\left[ e^{-\lambda u_{2}Z_{2}(t)\psi
(t)}\big |\mathbf{Z}(0)=\left( 0,N\right) \right] =e^{-ru_{2}\sqrt{\lambda }
} .  \label{inter3}
\end{equation}
Replacing \eqref{M1} by a stronger condition \eqref{FinFirst1} and assuming
that all the other assumptions hold we also have
\begin{equation}
\lim_{N,t\rightarrow \infty }\mathbf{E}\left[ e^{-\lambda u_{2}Z_{2}(t)\psi
(t)};Z_{1}(t)=0\big|\mathbf{Z}(0)=\left( 0,N\right) \right] =e^{-ru_{2}\sqrt{
\lambda }} .  \label{inter4}
\end{equation}
\end{theorem}

\noindent \textbf{Remark}. Recall that relation \eqref{VatSurviv}
implies that the population dies out whenever $N\sqrt{
1-G_{2}(t)}\rightarrow 0$, whereas the
limiting distribution of the number of particles is given by formula %
\eqref{Z12} whenever $N\sqrt{1-G_{2}(t)}\to r\in ( 0,\infty)$. In view of
the assumption $\beta<1$ the latter relation implies $R(t)/N\sim
(1-\beta)(N( 1-G_{2}(t)))^{-1} \to\infty$.

\subsection{Summary of the results obtained}

In this subsection we give a description of the splitting into the regions
in terms of the functions $\mathfrak{g}_{i}(N)$, $i=1,2,3$.

First we recall that if $\beta \in (0,1]$, and $t,N\rightarrow
\infty $ in such a way that\newline $t\sim \mathfrak{g}_{3}(N/r)$
or $t\gg \mathfrak{g}_{3}(N)$ (that is, if $t(N)$ belongs to one
of the final evolutionary stages of the process) then the
asymptotic behavior of the number of particles is given in
(\ref{Z12}) with $r=0$ in the second case. Under these conditions
the first type particles are absent in the limit. The same is true
for the second type particles if $r=0$, while the limiting
distribution of the second type particles for $r>0$ is discrete.
 We summarize this result in Table~0:

\begin{center}
\textbf{Table 0} case $\beta \in (0,1]$

\medskip

\begin{tabular}{|l|c|c|} \hline
Range of $t$ & $\simeq \mathfrak{g}_{3}$ & $\gg \mathfrak{g}_{3}$ \\
\hline
Theorems & $\text{f(\ref{Z12})}$ & $\text{f(\ref{Z12})}$ \\ \hline
\end{tabular}
\end{center}
where the first row shows the set of ranges of $t=t(N)$ under
consideration (with  $\mathfrak{g}_{i}$ for $\mathfrak{g}_{i}(N)$)
and "f(\ref{Z12})" in the second row indicates that the statement
for the corresponding time interval is given by formula
\eqref{Z12}.

Given below are three other tables which describe the evolution of the
two-type Bellman-Harris branching process in some remaining cases. Tables 1,
2 and 3 concern the cases $\beta \in (0,1/2)$, $\beta \in (1/2,2/3)$ and $%
\beta \in (2/3, 1)$, respectively.

\begin{center}
\textbf{Table 1} case $\beta \in (0,1/2)$

\medskip

\begin{tabular}{|l|c|c|c|c|c|} \hline Range of $t$ &
$o(\mathfrak{g}_{2})$ & $\simeq \mathfrak{g}_{2}$ & $ \in \left(
\mathfrak{g}_{2},\mathfrak{g}_{1}\right) $ & $\simeq
\mathfrak{g}_{1}$ & $\in \left( \mathfrak{g}_{1},\mathfrak{g}
_{3}\right) $ \\ \hline Theorems & $\text{T\ref{T_initial1}}$ &
$\text{ T\ref{T_initial2}}$ & $\text{T\ref{T_initial3}}$ & $\text{
T\ref{T_Zone3}}$ & $\text{T\ref{T_zone2}}$ \\
\hline
\end{tabular}
\end{center}


\begin{center}
\textbf{Table 2} case $\beta \in (1/2,2/3)$

\medskip

\begin{tabular}{|l|c|c|c|} \hline Range of $t$ &
$o(\mathfrak{g}_{1})$ & $\simeq \mathfrak{g}_{1}$ & $
\in \left( \mathfrak{g}_{1},\mathfrak{g}_{3}\right) $ \\
\hline
Theorems & $\text{T\ref{T_initial1}}$ & $\text{
T\ref{T_Betamore}}$ & $\text{T\ref{T_zone2}}$ \\ \hline
\end{tabular}
\end{center}

Here the symbols T\ref{T_initial1}, T\ref{T_initial2} etc. mean that the
result for the corresponding range of $t=t(N)$ is given in Theorem \ref%
{T_initial1}, Theorem \ref{T_initial2} etc. Thus, we have a complete
limiting picture for all fixed $\beta \in (0,2/3)\backslash\{1/2\}$ as $%
t=t(N)\rightarrow \infty $ (in some regions $t(N)$ should be a function
regularly varying at $\infty$). Moreover, provided that conditions (\ref%
{Redu}) or (\ref{Sup_cond}) hold Theorems \ref{T_initial1} and \ref%
{T_initial2} are true for $\beta =1/2$ as well, whereas the case $\beta =1/2$
is still open, otherwise.

\begin{center}
\textbf{Table 3} case $\beta \in (2/3,1)$

\medskip

\begin{tabular}{|l|c|c|c|} \hline Range of $t$ &
$o(\mathfrak{g}_{1})$ & $\simeq \mathfrak{g}_{1}$ & $ \in \left(
\mathfrak{g}_{1},\mathfrak{g}_{3}\right) $ \\ \hline Theorems &
$\text{T\ref{T_initial1}}$ & $\text{ T\ref{T_Betamore}}$ &
$\text{T\ref{T_zone2}}+?$ \\ \hline
\end{tabular}
\end{center}

Table 3 looks almost the same as Table 2. There are, however, two essential
differences. First, Table 2 concerns the case $\mathfrak{g}_{2}(N)\ll
\mathfrak{g}_{3}(N)$, i.e., $\beta<2(1-\beta)$, while Table 3 deals with the
case $\mathfrak{g}_{3}(N)\ll\mathfrak{g}_{2}(N)$, i.e., $\beta>2(1-\beta)$.
Second, in the case $\beta \in (2/3,1)$ and under certain conditions in the
case $\beta =2/3$ and $\mathfrak{g} _{2}(N)\leq t(N)\ll\mathfrak{g}_{3}(N)$
we have no results showing that the first type particles are absent within
the corresponding time interval (this fact is indicated by the sign ``+?''
in Table 3). We, however, believe that this is, indeed, the case.

\section{Auxiliary results\label{Sec2}}

We stipulate hereafter that, unless otherwise is stated, records like $%
a(t)=O(b(t))$, $a(t)=o(b(t))$, $a(t)\sim b(t)$ or $\lim a(t)=a$ are assumed
to hold, as $t\rightarrow \infty$.

Set $\mathbf{I}(t):=1_{\left\{ t\geq 0\right\} }\mathbf{I},$ where $%
1_{\{A\}} $ denotes the indicator of the event $A$, and $\mathbf{I}:\mathbf{%
= }\left( \delta _{ij}\right) _{i,j=1,2}$. We define the convolution $%
\mathbf{C }(t)=\mathbf{A}\ast \mathbf{B}(t)=\left( C_{ij}(t)\right)
_{i,j=1}^{2}$ of two matrices $\mathbf{A}(t)=\left( A_{ij}(t)\right)
_{i,j=1}^{2}$ and $\mathbf{B}(t)=\left( B_{ij}(t)\right) _{i,j=1}^{2}$ as
the matrix with elements
\begin{equation*}
C_{ij}(t)=\sum_{k=1}^{2}A_{ik}\ast B_{kj}(t).
\end{equation*}
Put $\mathbf{M}^{\ast 0}(t):=\mathbf{I}(t)$ and introduce the renewal matrix
\begin{equation*}
\mathbf{U}(t)=(U_{ij}(t))_{i,j=1}^{2}:=\sum_{k=0}^{\infty }\mathbf{M}^{\ast
k}(t)
\end{equation*}
with the agreement that $\mathbf{U}(t)$ is the $2\times 2$ zero matrix if $%
t<0$. Clearly,
\begin{equation}
\mathbf{U}(t)=\mathbf{I}(t)+\mathbf{M}\ast \mathbf{U}(t).  \label{MatRen}
\end{equation}
We also define the $L_{1}$ norm $\left\Vert \cdot \right\Vert $ of matrices
and vectors as the sum of the absolute values of all their components, i.e.,
\begin{equation*}
\left\Vert \mathbf{U}\left( t\right) \right\Vert =\sum_{i,j=1}^{2}U_{ij}(t).
\end{equation*}

The following statements concerning various asymptotic properties of the
renewal matrix $\mathbf{U}(t)$ have been established in \cite{VT13}.

\begin{lemma}
\label{L_asympRen} Under conditions \eqref{02} and \eqref{03}
\begin{equation*}
\mathbf{U}(t)\sim \Gamma _{\beta }\frac{t}{\mu _{2}(t)}\mathbf{D=}\Gamma
_{\beta }R(t)\mathbf{D}.
\end{equation*}
\end{lemma}

\begin{corollary}
\label{LRen} Suppose that conditions \eqref{02} and \eqref{03} hold. Then
there exists a constant $C\in \left( 0,\infty \right) $ such that, for all $%
t\geq 0,$
\begin{equation*}
\left\Vert \mathbf{U}\left( t\right) \right\Vert \leq C\left( R(t)+1\right) .
\end{equation*}
\end{corollary}

\begin{lemma}
\label{CShur02} Suppose that Hypothesis \textbf{A} holds. Then, for any
fixed $\Delta >0,$
\begin{equation}
\mathbf{U}(t+\Delta )-\mathbf{U}(t)\sim \Delta \frac{\beta \Gamma _{\beta }}{%
\mu _{2}(t)}\mathbf{D.}  \label{IncU1}
\end{equation}%
Moreover,
\begin{equation}
\mathbf{U}_{\mathbf{I}}(t):=\mathbf{U}\ast \mathbf{G}_{\mathbf{I}}(t)\sim
\mathbf{U}(t)  \label{Lu}
\end{equation}%
and
\begin{equation}
\mathbf{U}_{\mathbf{I}}(t+\Delta )-\mathbf{U}_{\mathbf{I}}(t)\sim \Delta
\frac{\beta \Gamma _{\beta }}{\mu _{2}(t)}\mathbf{D}.  \label{InU}
\end{equation}
\end{lemma}

\noindent {\textbf{Remark}. In Lemma~8 of \cite{VT13} relation
\eqref{InU} has only been stated for $\beta \in (0,1/2]$. A
perusal of the proof of that lemma reveals that \eqref{InU} holds
true for all $\beta \in (0,1]$.

\begin{lemma}
\label{CShur2} Assume that Hypothesis \textbf{A} holds. Then, for any
function $w(t)$ directly Riemann integrable on $[0,\infty )$
\begin{equation*}
\int_0^t w(t-u)\mathrm{d}\mathbf{U}(u) \sim \frac{\beta \Gamma _{\beta }}{
\mu _{2}(t)}\int_0^\infty w(u)\mathrm{d}u\mathbf{D},
\end{equation*}
and for any function $W(t)=t^{-\alpha }\ell_{W}(t)$, where $\ell_{W}(t)$ is
a function slowly varying at infinity,
\begin{eqnarray*}
\int_0^tW(t-u)\mathrm{d}\mathbf{U}(u) &\sim & {\frac{\beta \Gamma _{\beta
}\Gamma(1-\alpha)\Gamma(\beta)}{\Gamma(1-\alpha+\beta)}}{\frac{t^{1-\alpha
}\ell_{W}(t)}{\mu_2(t) }}\mathbf{D},\ \alpha \in \lbrack 0,1), \\
\int_0^t W(t-u)\mathrm{d}\mathbf{U}(u) &\sim &\frac{\beta \Gamma _{\beta }}{
\mu _{2}(t)}\int_0^t u^{-1}\ell_{W}(u)\mathrm{d}u\mathbf{D},\ \alpha =1.
\end{eqnarray*}
Here $\Gamma(\cdot)$ is the gamma function. In particular,
\begin{eqnarray*}
\int_0^t \left( 1-G_1(t-u\right) )\mathrm{d}\mathbf{U}(u) &\sim &\frac{\beta
\Gamma _{\beta }\mu_1}{\mu_2(t)}\mathbf{D}, \\
\int_0^t\left( 1-G_2(t-u\right) )\mathrm{d}\mathbf{U}(u) &\sim & \mathbf{D}.
\end{eqnarray*}
All the previous relations remain valid if we replace $\mathbf{U}(t)$ by $%
\mathbf{U}_{\mathbf{I}}(t)$.
\end{lemma}

The results of Lemma \ref{CShur2} concerning $\mathbf{U}$ were obtained in
Corollaries 7 and 9 in \cite{VT13}. The results concerning $\mathbf{U}_{
\mathbf{I}}$ are new and can be derived from Lemma \ref{CShur02} and Lemma 6
in \cite{VT13}.

Let
\begin{equation*}
{\mathbf{P}}(t):=\big(P_{ij}(t)\big)_{i,j=1,2}=\left(\mathbf{E}_{i}Z_{j}(t)
\right) _{i,j=1,2}
\end{equation*}
denote the mean matrix of the number of particles at time $t$. The
asymptotic behavior of $\mathbf{P}(t)$ is given by the following lemma.

\begin{lemma}
\label{L_firsrMom} $[($see $\cite{VT13})]$ Suppose that Hypothesis \textbf{A}
holds. Then
\begin{equation}
\mathbf{P}(t)\sim \left(
\begin{array}{cc}
D_{11}\frac{\textstyle\mu _{1}\beta \Gamma _{\beta }}{\textstyle\mu _{2}(t)+1%
} & D_{12}\medskip \\
D_{21}\frac{\textstyle\mu _{1}\beta \Gamma _{\beta }}{\textstyle\mu _{2}(t)+1%
} & D_{22}%
\end{array}
\right) =\mathbf{DJ}(t),  \label{AsymP}
\end{equation}
where
\begin{equation}
\mathbf{J}(t):=\left(
\begin{array}{cc}
\frac{\textstyle\mu _{1}\beta \Gamma _{\beta }}{\textstyle\mu _{2}(t)+1} &
0\medskip \\
0 & 1%
\end{array}
\right) .  \label{DefJ}
\end{equation}
\end{lemma}

Also, we mention that if
\begin{equation}
\mathbf{P}_{i}(t):=\left( P_{1i}(t),P_{2i}(t)\right) ^{\dagger },\ \mathbf{\
\delta }_{i}:=\left( \delta _{1i},\delta _{2i}\right) ^{\dagger },
\label{com0}
\end{equation}
then (see \cite{VT13}, formula (57) in Section 3)
\begin{equation}
\mathbf{P}(t)=\left( \mathbf{P}_{1}(t),\mathbf{P}_{2}(t)\right) =\mathbf{U}
\ast (\mathbf{I}-\mathbf{G}_{\mathbf{I}}(\cdot ))(t).  \label{pm2}
\end{equation}
In particular,
\begin{equation}
\mathbf{P}_{i}(t)=\mathbf{U}\ast \left( \mathbf{\delta }_{i}\otimes \left(
\mathbf{1}-\mathbf{G}(\cdot )\right) \right) (t),\quad i=1,2.  \label{com1}
\end{equation}

The section closes with two more technical lemmas of different flavor.

\begin{lemma}
\label{LRed} Let $H(t)$, $k_1(t)$ and $k_2(t)$ be nonnegative functions
defined on $[0,\infty)$ and such that $\lim_{t\rightarrow \infty }k_{i}(t)=0$%
, $i=1,2$, and $\lim_{t\to\infty}H(t)=\infty$. If
\begin{equation*}
\lim_{t\rightarrow \infty }H(t)\mathbf{Q}\left( t;1-\lambda_{1}k_{1}\left(
t\right) ,1-\lambda _{2}k_{2}(t)\right) =\mathbf{h}(\lambda _{1},\lambda
_{2}), \ \ \lambda\geq 0, \,\lambda_2\geq 0,
\end{equation*}
for a function $\mathbf{h}(\lambda_1,\lambda_2)=(h_1(\lambda_1,\lambda_2),
h_2(\lambda_1,\lambda_2))$ with continuous in both arguments components,
then
\begin{equation*}
\lim_{t\rightarrow \infty }H(t)\mathbf{Q}\left( t;e^{-\lambda _{1}k_3\left(
t\right) },e^{-\lambda _{2}k_4\left( t\right) }\right) =\mathbf{h}(\lambda
_{1},\lambda _{2})
\end{equation*}
for any functions $k_3(t)$ and $k_4(t)$ such that $\lim_{t\to%
\infty}k_3(t)/k_1(t)=1$ and $\lim_{t\to\infty}k_4(t)/k_2(t)=1$.
\end{lemma}

\textbf{Proof}. The result follows from the monotonicity of $Q_i\left(
t;s_{1},s_{2}\right)$, $i=1,2$, in $s_1$ and $s_2$, the inequality
\begin{equation}
1-x\leq e^{-x}\leq 1-x\left( 1-\frac{x}{2}\right) ,  \label{SIneq}
\end{equation}
and the continuity of $h_i (\lambda _{1},\lambda _{2})$, $i=1,2$, in both
arguments.

\begin{lemma}
\label{L_change2} Let $H(t)$, $k_{1}(t)$ and $k_{2}(t)$ be
nonnegative functions defined on $[0,\infty )$ and such that
$\lim_{t\rightarrow \infty }k_{i}(t)=0$, $i=1,2$, and
$\lim_{t\rightarrow \infty }H(t)=\infty $. If there exist
functions $k_i(t,\lambda_i), i=1,2$ such that
$k_i(t,\lambda_i)\sim \lambda_ik_i(t),\, t\to\infty$, for any
fixed $\lambda_i>0$, and
\[ \lim_{t\rightarrow \infty
}H(t)\mathbf{Q}\left( t;1-k_1(t,\lambda_1),1-k_2(t,\lambda_2)
\right) =\mathbf{h} (\lambda_1,\lambda_2), \ \lambda_1>0, \
\lambda_2>0,
\]
for a function $\mathbf{h}(\lambda _{1},\lambda
_{2})=(h_{1}(\lambda _{1},\lambda _{2}),h_{2}(\lambda _{1},\lambda
_{2}))$ with continuous in both arguments components, then
\begin{equation*}
\lim_{t\rightarrow \infty }H(t)\mathbf{Q}\left( t;e^{-\lambda
_{1}k_{1}(t)},e^{-\lambda _{2}k_{2}(t)}\right) =\mathbf{h}(\lambda
_{1},\lambda _{2}),\ \lambda _{1}>0,\ \lambda _{2}>0.
\end{equation*}%
If
\begin{equation*}
\lim_{t\rightarrow \infty }H(t)\mathbf{Q}\left( t;1,1-k_{2}(t,\lambda
_{2})\right) =\mathbf{h}(\lambda _{2}),\ \lambda _{2}>0,
\end{equation*}%
for a function $\mathbf{h}(\lambda _{2})=(h_{1}(\lambda _{2}),h_{2}(\lambda
_{2}))$ with continuous components, then
\begin{equation*}
\lim_{t\rightarrow \infty }H(t)\mathbf{Q}\left( t;1,e^{-\lambda
_{2}k_{2}(t)}\right) =\mathbf{h}(\lambda _{2}),\ \lambda _{2}>0.
\end{equation*}
\end{lemma}

\textbf{Proof}. We only prove the first part of the lemma. The
proof of the second part is analogous.

For any fixed positive $\lambda_1$, $\lambda_2$ and any
$\varepsilon \in \left(
0,\min\{1,\lambda_{1},\lambda_{2}\}/2\right)$ there exist
$T_i=T_i(\lambda_i), i=1,2 $ such that
\begin{equation}
\lambda_{i}k_{i}(t)> k_{i}(t,\lambda_i-\varepsilon),\quad
\lambda_{i}k_{i}(t)< k_{i}(t,\lambda_i+\varepsilon) \label{BLA}
\end{equation}
for $t\geq T_i$ which implies that both the upper and the lower
limits of $H(t)\mathbf{Q}( t;1-\lambda_1k_1(t),
1-\lambda_2k_2(t))$ are sandwiched between $\mathbf{h}
(\lambda_1+\varepsilon ,\lambda_2+\varepsilon )$ and $\mathbf{h}
(\lambda_1-\varepsilon ,\lambda_2-\varepsilon )$. Using the
continuity of the components of $\mathbf{h}$ in
$(\lambda_1,\lambda_2)$ yields $$\lim_{t\to\infty}H(t)\mathbf{Q}(
t;1-\lambda_1k_1(t), 1-\lambda_2k_2(t))=\mathbf{h}
(\lambda_1,\lambda_2),$$ and an appeal to Lemma \ref{LRed}
completes the proof of the first part.

\section{Proofs for the early evolutionary stages\label{Sec3}}

Recalling the definition $\mathbf{U}_{\mathbf{I}}(t)=\mathbf{U}\ast \mathbf{%
\ G }_{\mathbf{I}}(t)$ and viewing equation (\ref{EqvForQ}) as a
renewal-type equation with respect to $\mathbf{Q}(t;\mathbf{s})$ we obtain
the representation
\begin{equation}  \label{RenewalQ}
\mathbf{Q}(t;\mathbf{s})=\mathbf{U}\ast \left( \left( \mathbf{1}-\mathbf{s}
\right) \otimes \left( \mathbf{1}-\mathbf{G}(\cdot )\right) \right)
(t)-\int_{0}^{t}\mathrm{d}\mathbf{U}_{\mathbf{I}}(w)\mathbf{\Phi }(\mathbf{Q}
(t-w; \mathbf{s}))
\end{equation}
which is our main tool in this section.

\subsection{Proof of Theorem \protect\ref{T_initial1}}

Invoking Lemma \ref{LRed} which applies because $\mu_2(t)/N\to 0$ we
conclude that it suffices to verify that
\begin{equation}  \label{inter10}
\lim_{N,t\rightarrow \infty }N\mathbf{Q}\left( t;1-N^{-1}\lambda _{1}\mu
_{2}(t),1-N^{-1}\lambda _{2}\right) =\mathbf{D}\left( \lambda _{1}\mu
_{1}\beta \Gamma _{\beta },\lambda_{2}\right) ^{\dagger }.
\end{equation}

Put $\mathbf{s}:=\big(1-N^{-1}\lambda _{1}\mu _{2}(t), 1-N^{-1}\lambda _{2} %
\big)$. Our strategy is to show that
\begin{equation}  \label{inter12}
\lim_{N,t\rightarrow \infty }N\mathbf{U}\ast \left( \left( \mathbf{1}-
\mathbf{\ s}\right) \otimes \left( \mathbf{1}-\mathbf{G}(\cdot )\right)
\right)=\mathbf{D}\left( \lambda _{1}\mu _{1}\beta \Gamma _{\beta },\lambda
_{2}\right) ^{\dagger }
\end{equation}
and
\begin{equation}  \label{inter11}
N\int_0^t \mathrm{d}\mathbf{U}_{\mathbf{I}}(w)\mathbf{\Phi }(\mathbf{Q}(t-w;
\mathbf{\ s}))=o\left( \mathbf{1}\right)
\end{equation}
which entails \eqref{inter10} in view of \eqref{RenewalQ}.

Relation \eqref{inter12} is an immediate consequence of
\begin{eqnarray}
&&\mathbf{U}\ast \left( \left( \mathbf{1}-\mathbf{s}\right) \otimes \left(
\mathbf{1}-\mathbf{G}(\cdot )\right) \right) (t)=\mathbf{U}\ast \left(
\left( \mathbf{I}-\mathbf{G}_{I}(\cdot )\right) \right) (t)\left( \mathbf{1}%
- \mathbf{s}\right)  \notag \\
&&\qquad \qquad \qquad \qquad =\mathbf{P}(t)\left( \mathbf{1}-\mathbf{s}
\right) =\left( 1+o(1)\right) \mathbf{DJ}(t)\left( \mathbf{1}-\mathbf{s}
\right),  \label{AsMat}
\end{eqnarray}
where the second and the third equalities follow from \eqref{pm2} and Lemma %
\ref{L_firsrMom}, respectively.

Left with the proof of \eqref{inter11} we observe that \eqref{AsMat} ensures
the existence of $c\geq \left( \lambda _{1}+\lambda _{2}\right) $ such that
\begin{equation}
\mathbf{Q}(w;\mathbf{s})\leq \min \left\{ 1,c\frac{\mu _{2}(t)}{N\hat{\mu}
_{2}(w)}\right\} \mathbf{1}  \label{QQ}
\end{equation}
for all pairs $w\leq t$, where $\hat{\mu}_{2}(t):=\mu_{2}(t)+1$.
Further we recall that if $\Upsilon $ in \eqref{Sup_cond} is
finite then $\beta \in (0,0.5] $ whereas if $\Upsilon =\infty $
then $\beta \in [0.5,1]$. Writing
\begin{eqnarray*}
&&\int_0^t\mathrm{d}\mathbf{U}_{\mathbf{I}}(w)\mathbf{\Phi }(\mathbf{Q}(t-w;
\mathbf{s}))\leq \int_0^t \mathbf{U}_{\mathbf{I}}(\mathrm{d}w)\mathbf{\Phi }
\left( \min \left\{ 1,c\frac{\mu _{2}(t)}{N\hat{\mu}_{2}(t-w)}\right\}
\mathbf{1}\right)  \notag \\
&&\qquad \leq \ C\frac{\mu _{2}^{2}(t)}{N^{2}}\int_0^t\mathbf{U}_{\mathbf{I}
}(\mathrm{d}w)\frac{1}{\hat{\mu}_{2}^{2}(t-w)}\mathbf{1}  \notag \\
&&\qquad =O\left( \frac{\mathbf{1}\mu _{2}^{2}(t)}{N^{2}}\right) \times
\left\{
\begin{array}{lll}
\mu _{2}^{-1}(t), & \text{if} & \Upsilon <\infty , \\
\mu _{2}^{-1}(t)\int_{0}^{t}\frac{\textstyle{1}}{\textstyle\hat{\mu}
_{2}^{2}(u)}\mathrm{d}u, & \text{if} & \Upsilon =\infty ,\,\beta =0.5, \\
t\mu _{2}^{-3}(t), & \text{if} & \beta \in (0.5,1]%
\end{array}
\right .
\end{eqnarray*}
where the first inequality follows from \eqref{QQ}, the second is a
consequence of \eqref{DDEF} and \eqref{OcN}, and the third is justified by
Lemma \ref{CShur2}, and appealing to \eqref{Redu} in the case $\Upsilon
=\infty$ and $\beta =0.5$ and to \eqref{inter13} in the complementary cases
we arrive at \eqref{inter11}. The proof of Theorem \ref{T_initial1} is
complete.

\subsection{Proof of Corollary \protect\ref{C_initial1}}

Put $\mathbf{s}:=(1, 1-\lambda N^{-1})$. A similar argument as in the proof
of Theorem \ref{T_initial1} shows that it suffices to prove
\begin{equation}  \label{inter17}
\lim_{N,\,t\rightarrow \infty }N\mathbf{U}\ast \left( \left( \mathbf{1}-
\mathbf{\ s}\right) \otimes \left( \mathbf{1}-\mathbf{G}(\cdot )\right)
\right)=\mathbf{D}\left( 0,\lambda \right) ^{\dagger }
\end{equation}
and \eqref{inter11} (with the present $\mathbf{s}$).

Using \eqref{AsMat} (with the present $\mathbf{s}$) proves
\eqref{inter17}. As for \eqref{inter11}, mimicking the argument
given in the proof of Theorem \ref{T_initial1} leads to
$\mathbf{Q}(w;\mathbf{s})\leq \min \{1,cN^{-1}\}\mathbf{1}$,
$0\leq w\leq t$ and then to
\begin{eqnarray*}
\int_{0}^{t}\mathrm{d}\mathbf{U}_{\mathbf{I}}(w)\mathbf{\Phi }(\mathbf{Q}%
(t-w;\mathbf{\ s})) &\leq &\int_{0}^{t}\mathrm{d}\mathbf{U}_{\mathbf{I}}(w)%
\mathbf{\ \Phi }\left( \min \{1,cN^{-1}\}\mathbf{1}\right) \\
&=&O\left( N^{-2}\mathbf{U}_{\mathbf{I}}(t)\mathbf{1}\right) =o\left( N^{-1}%
\mathbf{1}\right) ,
\end{eqnarray*}%
where relation \eqref{Lu}, Lemma \ref{L_asympRen} and the assumption $%
R(t)/N\rightarrow 0$ have been utilized for the last equality. The proof of
Corollary \ref{C_initial1} is complete.

\subsection{Proof of Theorem \protect\ref{T_initial2}}

Put $\mathbf{s}:=\big(s, e^{-\lambda/N}\big)$. We intend to prove
\begin{equation}
\lim_{N,t\rightarrow \infty }N\mathbf{U}\ast \left( \left( \mathbf{1}-
\mathbf{s} \right) \otimes \left( \mathbf{1}-\mathbf{G}(\cdot )\right)
\right) (t)=\mathbf{D}\left(
\begin{array}{c}
r\left( 1-s\right) \mu _{1}\beta \Gamma _{\beta }\medskip \\
\lambda%
\end{array}
\right)  \label{inter15}
\end{equation}
and
\begin{equation}  \label{inter16}
\lim_{N,t\rightarrow \infty } N\int_0^t \mathrm{d}\mathbf{U}_{\mathbf{I}}(w)
\mathbf{\Phi }(\mathbf{Q}(t-w; \mathbf{s}))=r\mathbf{O}(s)
\end{equation}
which entails
\begin{equation*}
\lim_{N,\,t\rightarrow \infty }N\mathbf{Q}(t; s,e^{-\lambda/N})=\mathbf{D}
\left(
\begin{array}{c}
r\left( 1-s\right) \mu _{1}\beta \Gamma _{\beta } \\
\lambda%
\end{array}
\right) -r\mathbf{O}(s)
\end{equation*}
in view of \eqref{RenewalQ}.

Relation \eqref{inter15} is a consequence of \eqref{AsMat} and the
assumption $\mu_2(t)/N\to~r^{-1}$. Passing to the proof of \eqref{inter16}
we first note that \eqref{inter15} implies
\begin{equation}
Q_{i}(t;s,1)\leq Q_i(t;\mathbf{s})=O\left(1/N\right) =O\left( 1/\hat{\mu}
_{2}(t)\right)  \label{DD1}
\end{equation}
because $\mathbf{Q}(t-w;s,s_{2})$ is non-increasing in $s_2$. Since the
function $1/\hat{\mu}_{2}(t)$ is directly Riemann integrable on $[0,\infty)$
as a non-increasing Lebesgue integrable function (see \eqref{Sup_cond}), so
are $\Phi_i(\mathbf{Q}(t;s,1))$, $i=1,2$ because these are nonnegative
bounded and continuous functions satisfying
\begin{equation*}
\Phi_i(\mathbf{Q}(t;s,1))\leq C(Q_1^2(t;s,1)+Q_2^2(t;s,1))\leq C_1/\hat{\mu}
_{2}(t),
\end{equation*}
where the first inequality is a consequence of \eqref{DDEF} and \eqref{OcN},
and the second follows from \eqref{DD1}. With this at hand an application of
Lemma \ref{CShur2} yields
\begin{eqnarray}
&&\lim_{N,t\rightarrow \infty }N\int_0^t \mathrm{d}\mathbf{U}_{\mathbf{I}%
}(w) \mathbf{\Phi }(\mathbf{Q}(t-w;s,1))  \notag \\
&&\qquad \qquad =r\beta \Gamma _{\beta }\int_{0}^{\infty }\mathbf{D\Phi }(
\mathbf{Q}(w;s,1))\mathrm{d}w=r\mathbf{O}(s)<\infty .  \label{Two}
\end{eqnarray}
It remains to check that \eqref{Two} implies \eqref{inter16}. To this end,
write
\begin{equation*}
\mathbf{Q}(w; s,e^{-\lambda /N})-\mathbf{Q}(w;s,1)\leq \lambda N^{-1}\mathbf{%
\ P}_2(w)\leq CN^{-1}\mathbf{1},
\end{equation*}
where the last inequality follows from \eqref{AsymP}, and then
\begin{eqnarray*}
0 &\leq &\mathbf{\Phi }(\mathbf{Q}(w;s,e^{-\lambda /N}))-\mathbf{\Phi }(
\mathbf{Q}(w;s,1)) \\
&\leq& C\left( \mathbf{Q}(w;s,e^{-\lambda /N})-\mathbf{Q}(w;s,1)\right)
\left\Vert \mathbf{Q}(w;s,e^{-\lambda /N})\right\Vert \\
&\leq& \frac{C}{N}\left\Vert \mathbf{Q}(w;s,e^{-\lambda /N})\right\Vert
\mathbf{1}\leq \frac{C}{N\hat{\mu}_{2}(w)}\mathbf{1}
\end{eqnarray*}
to infer
\begin{eqnarray*}
0&\leq& N\int_0^t \mathrm{d}\mathbf{U}_{\mathbf{I}}(w)\left[ \mathbf{\Phi }(
\mathbf{Q}(t-w;s,e^{-\lambda /N}))-\mathbf{\Phi }(\mathbf{Q}(t-w;s,1))\right]
\\
&\leq& C\int_0^t \frac{\mathrm{d}\mathbf{U}_{\mathbf{I}}(w)}{\hat{\mu}
_2\left( t-w\right) }.
\end{eqnarray*}
Using Lemma \ref{CShur2} and the conditions of the theorem we conclude
\begin{equation*}
\int_0^t\frac{\mathrm{d}\mathbf{U}_{\mathbf{I}}(w)}{\hat{\mu}_2\left(
t-w\right) }=O\left( t\mu_2^{-2}(t)\right) =O\left( R(t)/\mu _{2}(t)\right)
=o(1)
\end{equation*}
which completes the proof of Theorem~\ref{T_initial2}.

\subsection{Proof of Theorem \protect\ref{T_initial3}}

The second equality in \eqref{inter14} has already been verified
in the proof of Corollary~\ref{C_initial1} under the sole
assumption $R(t)/N\to 0$. With this at hand the first equality in
\eqref{inter14}, equivalently,
\begin{equation*}
\lim_{N,t\rightarrow \infty }N\mathbf{Q}(t;0,
1-\lambda/N)=\exp(-D_{22}\lambda)
\end{equation*}
follows from the estimate
\begin{eqnarray*}
\mathbf{0} &\leq &N\bigg(\mathbf{Q}(t;0,1-\lambda/N)-\mathbf{Q}
(t;1,1-\lambda/N)\bigg)  \notag \\
&\leq &N \mathbf{P}_{1}\left( t\right)=O(N/\mu_2(t))=o(1),
\end{eqnarray*}
where Markov's inequality has been used for the second inequality and Lemma %
\ref{L_firsrMom} and the assumption $\mu_2(t)/N\to\infty$ for the first and
the second equalities, respectively. The proof of Theorem \ref{T_initial3}
is complete.

\section{Proofs for the intermediate evolutionary stages\label{Sec4}}

The intermediate evolutionary stages exhibit the most interesting
and exotic behavior. Our main technical tool here is the
Contraction Principle.

\subsection{Proof of Theorem \protect\ref{T_Zone3}}

Lemma \ref{L_uniq} given below is an important ingredient of the proof of
Theorem~\ref{T_Zone3}.

\begin{lemma}
\label{L_uniq} For each $\beta \in (0,1]$ there exists $\Lambda >0$ such
that in the domain
\begin{equation*}
\mathcal{K}=\mathcal{K}\left( \beta ,\Lambda \right) :=\left\{ 0<\lambda
\leq \Lambda \right\} \subset \mathbb{R}
\end{equation*}
equation \eqref{UniNew1} has a unique solution in the class of
vector-functions with non-negative continuous components.
\end{lemma}

\noindent \textbf{Proof.} We set
\begin{equation*}
\mathbf{\Theta }\left( \lambda \right) =\left( \Theta _{1}(\lambda ),\Theta
_{2}(\lambda )\right) ^{\dagger }:=\lambda ^{-\beta }\mathbf{\Omega }\left(
\lambda ^{\beta }\right) =\left( \lambda ^{-\beta }\Omega _{1}(\lambda
^{\beta }),\lambda ^{-\beta }\Omega _{2}(\lambda ^{\beta })\right) ^{\dagger}
\end{equation*}
and consider the following system which is equivalent to (\ref{UniNew1})
\begin{equation}
\mathbf{\Theta }\left( \lambda \right) =\mathbf{D}\left( 0,1\right)
^{\dagger }-\Gamma _{\beta }\lambda ^{\beta }\int_{0}^{1}\mathbf{DN}\left(
\mathbf{\Theta }\left( \lambda \left( 1-y\right) \right) \right) \mathrm{d}
y^\beta.  \label{UNiNew33}
\end{equation}

If a desired solution of (\ref{UniNew1}) exists, then $\mathbf{\Theta }
\left( \lambda \right)$ is a fixed point of the respective mapping. Thus, it
is natural to approximate this fixed point by a sequence of iterates. To
this end, starting with $\mathbf{\Theta }^{(0)}\left( \lambda \right) :=
\mathbf{D}\left( 0,1\right) ^{\dagger }$ we define
\begin{eqnarray*}
\mathbf{\Theta }^{(n+1)}\left( \lambda \right)& :=&\mathbf{D}\left(
0,1\right) ^{\dagger }-\Gamma _{\beta }\lambda ^{\beta }\int_0^1\mathbf{DN}
\left( \mathbf{\Theta }^{(n)}\left( \lambda \left( 1-y\right) \right)
\right) \mathrm{d}y^\beta \\
&=&\left( \Theta _{1}^{(n+1)}(\lambda ),\Theta _{2}^{(n+1)}(\lambda )\right)
^{\dagger }
\end{eqnarray*}
for $n=0,1,\ldots$

The set $\mathcal{C}_+[0,\Lambda]$ of continuous functions $h: [0,\Lambda]
\mapsto [0,\infty)\times[0,\infty)$ equipped with the metric $%
\rho(h_1,h_2)=\sup_{\lambda \in [0,\Lambda]}\left\Vert h_1(\lambda
)-h_2(\lambda )\right\Vert$ is a complete metric space. Invoking the
Contraction Principle we conclude that it suffices to prove
\begin{equation}  \label{inter18}
\mathbf{\Theta }^{(n)}\left( \lambda \right)\in \mathcal{C}_+[0,\Lambda]
\end{equation}
for $n=0,1,\ldots$ and
\begin{equation}  \label{inter19}
\rho(\mathbf{\Theta }^{(n+1)}, \mathbf{\Theta }^{(n)})\leq \kappa \rho(
\mathbf{\Theta }^{(n)}, \mathbf{\Theta }^{(n-1)})
\end{equation}
for $n=1,2,\ldots$ and appropriate $\kappa\in (0,1)$.

In view of
\begin{equation*}
\lambda ^{\beta }\int_{0}^{1}\mathbf{DN}\left( \mathbf{\Theta }^{(0)}\left(
\lambda _{y}\right) \right) \mathrm{d}y^{\beta }\leq C\lambda ^{\beta }
\mathbf{D}\int_{0}^{1}\mathrm{d}y^{\beta }\leq C_{1}\lambda ^{\beta } %
\mathbf{D}\left( 0,1\right) ^{\dagger }.
\end{equation*}
where $\lambda _{y}:=\lambda \left( 1-y\right)$, there exists a sufficiently
small $\Lambda \in \left( 0,1\right) $ such that
\begin{equation*}
\mathbf{\Theta }^{(1)}\left( \lambda \right) =\mathbf{D}\left( 0,1\right)
^{\dagger }-\Gamma _{\beta }\lambda ^{\beta }\int_{0}^{1}\mathbf{DN}\left(
\mathbf{\Theta }^{(0)}\left( \lambda _{y}\right) \right) \mathrm{d}y^{\beta
}\geq \mathbf{0}
\end{equation*}
for all $\lambda \in [0,\Lambda]$. Assume now that $\mathbf{\Theta }
^{(n)}\left( \lambda \right) \geq \mathbf{0}$ for $\lambda\in [0,\Lambda]$.
Using the estimate
\begin{equation*}
\mathbf{\Theta }^{(n)}\left( \lambda \right) \leq \mathbf{D}\left(
0,1\right) ^{\dagger }, \ \ \lambda\in [0,\Lambda]
\end{equation*}
and arguing as above we conclude $\mathbf{\Theta }^{(n+1)}\left( \mathbf{\
\lambda }\right)\geq \mathbf{0}$ for all $\lambda\in [0,\Lambda]$.
Continuity of components of $\mathbf{\Theta }^{(n)}\left( \mathbf{\lambda }
\right)$, $n=0,1,\ldots$ is obvious, and \eqref{inter18} follows.

Further we have
\begin{eqnarray*}
\left\Vert \mathbf{\Theta }^{(n+1)}(\lambda )-\mathbf{\Theta }^{(n)}(\lambda
)\right\Vert &\leq &C\lambda ^{\beta }\int_{0}^{1}\left\Vert \mathbf{N}
\left( \mathbf{\Theta }^{(n)}(\lambda _{y})\right) -\mathbf{N} \left(
\mathbf{\Theta }^{(n-1)}(\lambda _{y})\right) \right\Vert \mathrm{d}y^{\beta
} \\
&\leq &C_{1}\lambda ^{\beta }\int_{0}^{1}\left\Vert \mathbf{\Theta }^{(n)}
(\lambda _{y})-\mathbf{\Theta }^{(n-1)}(\lambda _{y})\right\Vert \mathrm{d}
y^{\beta } \\
&\leq& C_{1}\lambda ^{\beta }\sup_{\eta \in \mathcal{K}}\left\Vert \mathbf{\
\Theta }^{(n)}(\eta )-\mathbf{\Theta }^{(n-1)}(\eta )\right\Vert
\end{eqnarray*}
for all $\lambda\in [0,\Lambda]$, having utilized \eqref{NorN} and $\mathbf{%
0 }\leq \mathbf{\Theta }^{(n)}\left( \mathbf{\lambda }\right) \leq \mathbf{D}%
\left(0,1\right) ^{\dagger }$, $\lambda\in [0,\Lambda]$, $n=0,1,\ldots$ for
the second inequality, and thereupon \eqref{inter19} with $%
\kappa=C_{1}\Lambda ^{\beta }$ for $\Lambda>0$ such that $\kappa<1$. The
proof of Lemma \ref{L_uniq} is complete.

\begin{lemma}
\label{L_z3} Suppose that Hypothesis \textbf{A} holds, and
\begin{equation*}
\lim_{N,\,t\rightarrow \infty }R(t)N^{-1}=1/r\in \left( 0,\infty \right) .
\end{equation*}%
Then, for any $\lambda >0$,
\begin{equation}
\lim_{N,t\rightarrow \infty }N\mathbf{Q}\left( t;1,e^{-\lambda r/N}\right) =r%
\mathbf{\Omega }(\lambda ),  \label{cpp}
\end{equation}%
where $\mathbf{\Omega }(\lambda )\geq 0$ solves equation
\eqref{UniNew1}. Moreover, equation \eqref{UniNew1} has a unique
analytic solution (being a complex-valued vector-function) in the
half-plane $\lambda\in \mathbb{C}$, $\Re \lambda>0$.
\end{lemma}

\noindent \textbf{Proof.} For $\lambda >0$, set
\begin{equation}
\mathbf{K}(t,w;\lambda)=\left( K_{1}(t,w;\lambda ),K_{2}(t,w;\lambda
)\right) :=R(t/\lambda )\mathbf{Q}\left( w;1,1-\frac{1}{R(t/\lambda )}\right)
\label{DefK_i}
\end{equation}%
and
\begin{equation*}
\mathbf{K}(t;\lambda )=\left( K_{1}(t;\lambda ),K_{2}(t;\lambda )\right) :=%
\mathbf{K}(t,t;\lambda )=R(t/\lambda )\mathbf{Q}\left( t;1,1-\frac{1}{%
R(t/\lambda )}\right) .
\end{equation*}
It suffices to prove the existence of a vector-function
$\mathbf{K}$ which, in a domain $\mathcal{K}_1\subseteq
\mathcal{K}$ to be specified later, satisfies
\begin{equation}
\mathbf{K}(\lambda )=\mathbf{D}\left( 0,1\right) ^{\dagger }-\Gamma _{\beta
}\lambda ^{\beta }\int_{0}^{1}\mathbf{DN}\left( \mathbf{K}(\lambda
(1-w))\right) \mathrm{d}w^{\beta }\   \label{inter20}
\end{equation}%
and
\begin{equation}
\lim_{t\rightarrow \infty }\mathbf{K}(t;\lambda )=\mathbf{K}(\lambda )\ .
\label{inter21}
\end{equation}%
Indeed, by Lemma \ref{L_uniq} equation \eqref{UNiNew33} then has a
unique solution $\mathbf{\Theta }(\lambda )=\mathbf{K}(\lambda )$
for $\lambda\in\mathcal{K}_1$ which implies that equation
\eqref{UniNew1} has a unique solution $\mathbf{\Omega }\left(
\lambda \right) =\lambda \mathbf{\Theta }(\lambda ^{1/\beta })$
for $\lambda\in\mathcal{K}_1$. Furthermore, \eqref{inter21}
entails
\begin{equation*}
\lambda \lim_{t\rightarrow \infty }R(t/\lambda ^{1/\beta
})\mathbf{Q}\left( t;1,1-\frac{1}{R(t/\lambda ^{1/\beta })}\right)
=\mathbf{\Omega }\left( \lambda \right), \ \
\lambda\in\mathcal{K}_1
\end{equation*}%
and thereupon
\begin{equation*}
\lim_{t\rightarrow \infty }N\mathbf{Q}\left( t;1,1-{r\lambda\over
N}\right) =r\mathbf{\Omega }\left( \lambda \right), \ \ \lambda\in
\mathcal{K}_1
\end{equation*}%
which implies \eqref{cpp} for $\lambda\in\mathcal{K}_{1}$
in view of Lemma \ref{L_change2}. Hence we get%
\begin{equation}
\lim_{N,\,t\rightarrow \infty }F_{i}^{N}\left( t;1,e^{-\lambda
r/N}\right) =e^{-r\Omega _{i}\left( \lambda \right) }, \ \ \lambda
\in \mathcal{K}_1, \ i=1,2. \label{LALA}
\end{equation}%

Since $F_i^N\left( t;1,e^{-\lambda r/N}\right)$, $i=1,2$ are the
Laplace transforms of nonnegative random variables and the limits
exist for $\lambda \in \mathcal{K}_{1},$ it follows by the
uniqueness theorem for Laplace transforms that the limits at the
left-hand side of \eqref{LALA} exist for all $\lambda>0$.
Moreover, the limits are Laplace transforms and, therefore, there
exists a function $\mathbf{\Omega }^{\ast }(\lambda
)=(\Omega^\ast_1(\lambda), \Omega^\ast_2(\lambda))$, $\lambda
>0$ such that
\begin{equation*}
\lim_{N,t\rightarrow \infty }N\mathbf{Q}\left( t;1,e^{-\lambda
r/N}\right) =r\mathbf{\Omega }^{\ast }(\lambda ), \ \ \lambda
>0.
\end{equation*}
Also, for $i=1,2$, the function $\Omega_i^\ast\left( \lambda
\right)$ is analytic for $\Re\lambda
>0$ (under an appropriate choice of the branches in the case of
singularities) and such that $\Omega_i^\ast \left( \lambda
\right)=\Omega_i(\lambda )$ for $\lambda \in \mathcal{K}_1$. Since
$\mathbf{\Omega }^\ast(\lambda )$ solves \eqref{UniNew1} for
$\lambda \in \mathcal{K}_1$, it follows from the uniqueness
theorem for analytic functions that $\mathbf{\Omega }^\ast(\lambda
)$ solves \eqref{UniNew1} for all $\lambda\in \mathbb{C}$ with
$\Re\lambda
>0$ and, therefore, coincides with $\mathbf{\Omega }\left( \lambda
\right) $ as desired.

Thus, we concentrate on proving (\ref{inter21}) and check that
\begin{equation}
\lim_{t\rightarrow \infty }\sup_{\lambda \in \mathcal{K}_{1}}\left\Vert
\mathbf{K}(t;\lambda )-\mathbf{\Theta }\left( \lambda \right) \right\Vert =0
\label{inter22}
\end{equation}%
for $\mathbf{\Theta }$ defined in the proof of Lemma \ref{L_uniq}. To this
end, we use \eqref{RenewalQ} to obtain
\begin{eqnarray}
\mathbf{K}(t;\lambda ) &=&R(t/\lambda )\mathbf{U}\ast \left( \left( \mathbf{1%
}-\mathbf{s}\right) \otimes \left( \mathbf{1}-\mathbf{G}(\cdot )\right)
\right) (t)  \notag \\
&&-R(t/\lambda )\int_{0}^{1}\mathrm{d}\mathbf{U}_{\mathbf{I}}(ty)\mathbf{%
\Phi }(\mathbf{Q}(t_{y};\mathbf{s})),\ \ \lambda>0,
\label{inter23}
\end{eqnarray}%
where $\mathbf{s}=(s_{1},s_{2}):=\big(1,1-1/R(t/\lambda )\big)$ and $%
t_{y}:=t(1-y).$

Using \eqref{AsMat} and Lemma \ref{L_firsrMom} gives
\begin{eqnarray}
R(t/\lambda )\mathbf{U}\ast \left( \left( \mathbf{1}-\mathbf{s}\right)
\otimes \left( \mathbf{1}-\mathbf{G}(\cdot )\right) \right) (t) &=&(1+o(1))%
R(t/\lambda )\mathbf{DJ}(t)\left( \mathbf{1}-\mathbf{s}%
\right)  \notag \\
&=&\mathbf{D}\left( 0,1\right) ^{\dagger }(1+\delta _{0}(t;\lambda )),
\label{MM1}
\end{eqnarray}%
where $\lim_{t\rightarrow \infty }\delta _{0}(t;\lambda )=0$ uniformly in $%
\lambda >0$. We have
\begin{eqnarray*}
\Delta _{1}(t;\varepsilon ,\lambda )&:= &R(t/\lambda )\left\Vert
\int_{1-\varepsilon }^{1}\mathrm{d}\mathbf{U}_{\mathbf{I}}(ty)\mathbf{\Phi }(%
\mathbf{Q}(t_{y};\mathbf{s}))\right\Vert \\
&\leq &CR(t/\lambda )\left\Vert \int_{1-\varepsilon }^{1}\mathrm{d}\mathbf{U}%
_{\mathbf{I}}(ty)\mathbf{\Phi }\left( \frac{1}{R(t/\lambda )}\mathbf{1}%
\right) \right\Vert \\
&\leq &\frac{C}{R(t/\lambda )}\left\Vert \mathbf{U}_{\mathbf{I}}(t)-\mathbf{U%
}_{\mathbf{I}}(t\left( 1-\varepsilon \right) )\right\Vert \\
&\sim &C\lambda ^{\beta }\big(1-(1-\varepsilon )^{\beta }\big)
\end{eqnarray*}%
for any $\varepsilon \in \left( 0,1\right)$ and all $\lambda %
>0$ which proves
\begin{equation*}
\underset{\varepsilon \downarrow 0}{\lim }\,\underset{t\rightarrow \infty }{%
\lim \sup }\,\sup_{\lambda \in \mathcal{K}}\Delta _{1}(t;\varepsilon
,\lambda )=0.
\end{equation*}%
While the fourth line of the last displayed formula follows from %
\eqref{012aa} and Lemma \ref{L_asympRen}, the second is a consequence of
\begin{equation}
R(t/\lambda )\mathbf{Q}(w;\mathbf{s})\leq R(t/\lambda )\left( 1-s_{2}\right)
\mathbf{P}_{2}(w)=\mathbf{P}_{2}(w)\leq C\mathbf{1},\ \ w\in \lbrack 0,t]
\label{EstK_i}
\end{equation}%
which, in its turn, is justified by Lemma \ref{L_firsrMom}.

We have, for all $\lambda \in \mathcal{K}$,
\begin{eqnarray*}
&&R(t/\lambda )\int_{0}^{1-\varepsilon }\mathrm{d}\mathbf{U}_{\mathbf{I}}(ty)%
\mathbf{\Phi }\left( \mathbf{Q}(t_{y};\mathbf{s})\right) \\
&&\qquad=R(t/\lambda )\int_{0}^{1-\varepsilon }\mathrm{d}\mathbf{U}_{\mathbf{I}%
}(ty)\mathbf{\ \Phi }\left( \frac{1}{R(t/\lambda
)}\mathbf{K}\left(
t_{y};\lambda _{y}\right) \right) \\
&&\qquad=R(t/\lambda )(1+\delta _{1}\left( t;\lambda \right)
)\int_{0}^{1-\varepsilon }\mathrm{d}\mathbf{U}_{\mathbf{I}}(ty)\mathbf{N}%
\left( \frac{1}{R(t/\lambda )}\mathbf{K}\left( t_{y};\lambda _{y}\right)
\right) \\
&&\qquad=\frac{R(t)}{R(t/\lambda )}(1+\delta _{1}\left( t;\lambda
\right)
)\int_{0}^{1-\varepsilon }\frac{\mathrm{d}\mathbf{U}_{\mathbf{I}}(ty)}{R(t)}%
\mathbf{N}\left( \mathbf{K}\left( t_{y};\lambda _{y}\right) \right) \\
&&\qquad=\left( \lambda ^{\beta }+\delta _{2}\left( t;\lambda
\right) \right)
\int_{0}^{1-\varepsilon }\frac{\mathrm{d}\mathbf{U}_{\mathbf{I}}(ty)}{R(t)}%
\mathbf{N}\left( \mathbf{K}\left( t_{y};\lambda _{y}\right) \right)
\end{eqnarray*}%
having utilized \eqref{DDEF}, \eqref{EstK_i}, \eqref{Lu} and Lemma \ref%
{L_asympRen} for the third line. Furthermore,
\begin{equation*}
\lim_{t\rightarrow \infty }\sup_{\lambda \in \mathcal{K}}|\delta
_{1}(t;\lambda )|=\lim_{t\rightarrow \infty }\sup_{\lambda \in \mathcal{K}%
}|\delta _{2}(t;\lambda )|=0
\end{equation*}%
because $\lim_{t\rightarrow \infty }\sup_{\lambda \in \mathcal{K}%
}|R(t)/R(t/\lambda )-\lambda ^{\beta }|=0$ by Theorem 1.5.2 in \cite{BGT}.
Combining pieces together gives
\begin{eqnarray}
\mathbf{K}(t;\lambda ) &=&\mathbf{D}\left( 0,1\right) ^{\dagger }(1+\Delta
_{2}(t;\varepsilon ,\lambda ))  \notag \\
&&-\left( \lambda ^{\beta }+\delta _{2}\left( t;\lambda \right) \right)
\int_{0}^{1-\varepsilon }\frac{\mathrm{d}\mathbf{U}_{\mathbf{I}}(ty)}{R(t)}%
\mathbf{N}\left( \mathbf{K}\left( t_{y};\lambda _{y}\right) \right) ,
\label{ApprK2}
\end{eqnarray}%
where
\begin{equation*}
\lim_{\varepsilon \downarrow 0}\mathop{\lim\sup}\limits_{t\rightarrow \infty
}\sup_{\lambda \in \mathcal{K}}\left\vert \Delta _{2}\left( t;\varepsilon
,\lambda \right) \right\vert =0.
\end{equation*}

According to the proof of Lemma \ref{L_uniq} each component of the
vector-function $y\mapsto \mathbf{N}\left( \mathbf{\Theta }\left(
\lambda(1-y)\right)\right) $, $y\in [0,1]$, $\lambda\in \mathcal{K}$ is
nonnegative, bounded and continuous. This in combination with Lemma \ref%
{L_asympRen} gives
\begin{equation*}
\lambda ^\beta \Gamma_\beta \int_0^1 \mathbf{DN}\left( \mathbf{\Theta }
\left( \lambda _{y}\right) \right) \mathrm{d}y^\beta=\lim_{\varepsilon
\downarrow 0}\lambda ^{\beta }\int_0^{1-\varepsilon }\frac{\mathrm{d}\mathbf{%
\ U}_{\mathbf{I}}(ty)}{R(t)}\mathbf{N}\left( \mathbf{\Theta }\left( \lambda
_{y}\right)\right)
\end{equation*}
allowing to rewrite equation (\ref{UNiNew33}) for $\lambda\in\mathcal{K}$ as
\begin{equation}
\mathbf{\Theta }\left( \lambda \right) =\mathbf{D}\left( 0,1\right)
^{\dagger }\left( 1+\Delta _{3}\left( t;\varepsilon ,\lambda \right) \right)
-\lambda ^{\beta }\int_{0}^{1-\varepsilon }\frac{\mathrm{d}\mathbf{U}_{
\mathbf{I}}(ty)}{R(t)}\mathbf{N}\left( \mathbf{\Theta }\left( \lambda
_{y}\right) \right),  \label{ApprPsi222}
\end{equation}
where
\begin{equation*}
\lim_{\varepsilon \downarrow 0}\mathop{\lim\sup}\limits_{t\rightarrow \infty
}\sup_{\lambda \in \mathcal{K}}\left\vert \Delta _{3}\left( t;\varepsilon
,\lambda \right) \right\vert =0.
\end{equation*}

Using \eqref{ApprK2} and \eqref{ApprPsi222} we infer
\begin{eqnarray*}
&&\left\Vert \mathbf{K}(t;\lambda )-\mathbf{\Theta }\left( \lambda \right)
\right\Vert \\
&&\qquad \leq \lambda ^{\beta }\left\Vert \int_{0}^{1-\varepsilon }\frac{
\mathrm{d}\mathbf{U}_{\mathbf{I}}(ty)}{R(t)}(\mathbf{N}\left( \mathbf{K}
(t_{y};\lambda _{y})\right) -\mathbf{N}\left( \mathbf{\Theta }\left( \lambda
_{y}\right) \right) )\right\Vert +\Delta _{4}(t;\varepsilon ) \\
&&\qquad \leq C_{1}\lambda ^{\beta }\int_{0}^{1-\varepsilon }\left\Vert
\frac{\mathrm{d}\mathbf{U}_{\mathbf{I}}(ty)}{R(t)}\right\Vert \left\Vert
\mathbf{K}(t_{y};\lambda _{y})-\mathbf{\Theta }(\lambda _{y})\right\Vert
+\Delta _{4}(t;\varepsilon ) \\
&&\qquad \leq C_{2}\lambda ^{\beta }\sup_{y\geq t\varepsilon }\sup_{\theta
\in \mathcal{K}}\left\Vert \mathbf{K}(y;\theta ))-\mathbf{\Theta }(\theta
)\right\Vert \int_{0}^{1-\varepsilon }\left\Vert \frac{\mathrm{d}\mathbf{U}%
_{ \mathbf{I}}(ty)}{R(t)}\right\Vert +\Delta _{4}(t;\varepsilon ) \\
&&\qquad \leq C_{3}\lambda ^{\beta }\sup_{w\geq t\varepsilon }\sup_{\theta
\in \mathcal{K}}\left\Vert \mathbf{K}(w;\theta )-\mathbf{\Theta }(\theta
)\right\Vert +\Delta _{4}(t,\varepsilon )
\end{eqnarray*}
for $\lambda \in \mathcal{K}$, where the inequality
\begin{equation*}
\left\Vert \int_{0}^{1-\varepsilon }\frac{\mathrm{d}\mathbf{U}_{\mathbf{I}
}(tw)}{R(t)}\mathbf{N}\left( \mathbf{K}\left( t_{y};\lambda _{y}\right)
\right) \right\Vert \leq C\left\Vert \int_{0}^{1-\varepsilon }\frac{\mathrm{%
d }\mathbf{U}_{\mathbf{I}}(ty)}{R(t)}\right\Vert \leq C_{1}
\end{equation*}
that follows from Lemma \ref{L_asympRen} and (\ref{Lu}) has been utilized
for the second line, the third line following from (\ref{NorN}) and the last
from Lemma \ref{L_asympRen}. Observe that the constants $C_{i}$, $i=1,2,3$
do not depend on $\lambda \in \mathcal{K}$ and
\begin{equation*}
\Delta _{4}(t,\varepsilon )\leq C\sup_{\lambda \in \mathcal{K}}\left(
\left\vert \Delta _{2}\left( t;\varepsilon ,\lambda \right) \right\vert
+\left\vert \Delta _{3}\left( t;\varepsilon ,\lambda \right) \right\vert
+\left\vert \delta _{2}\left( t;\lambda \right) \right\vert \right) .
\end{equation*}
By shrinking, if needed, $\mathcal{K}$ (which gives $\mathcal{K}_1$) we can
and do assume that $C_{3}\sup\limits_{\lambda \in \mathcal{K}_1}\lambda
^{\beta }=C_{3}\Lambda ^{\beta }=\kappa <1$. Setting
\begin{equation*}
S(t):=\sup_{w\geq t}\sup_{\theta \in \mathcal{K}_1}\left\Vert \mathbf{K}
(w;\theta )-\mathbf{\Theta }(\theta )\right\Vert
\end{equation*}
and letting first $t\rightarrow \infty $ and then $\varepsilon \downarrow 0$
in the inequality
\begin{equation*}
S(t)\leq \kappa S(t\varepsilon )+\sup_{s\geq t}\Delta _{4}(s;\varepsilon )
\end{equation*}
we arrive at \eqref{inter22}. The proof of Lemma \ref{L_z3} is
herewith finished.

Turning to the proof of Theorem \ref{T_Zone3} it only takes to observe that %
\eqref{RconsGlobal} is an immediate consequence of Lemma \ref{L_z3} and
relation \eqref{Blarge}, whereas \eqref{Cr} follows from
\begin{equation*}
\mathbf{0}\leq N\big(\mathbf{Q}(t;0,e^{-\lambda r/N })-\mathbf{Q}
(t;1,e^{-\lambda r/N})\big)\leq N\mathbf{P}_{1}(t)=O(N/\mu_2(t))\mathbf{1}
=o(1)\mathbf{1},
\end{equation*}
where Lemma \ref{L_firsrMom} has been used for the penultimate equality, and
the assumption $\mu_2(t)/N\to\infty$ for the last. The proof of Theorem \ref%
{T_Zone3} is complete.

\subsection{Proof of Theorem \protect\ref{T_Betamore}}

The proof of Theorem \ref{T_Betamore} rests on the following lemma.

\begin{lemma}
\label{L_fullShur} For each $\beta \in (1/2,1]$ there exists $\Lambda >0$
such that in the domain
\begin{equation*}
\mathcal{K}=\mathcal{K}(\beta ,\Lambda ,1):=\left\{ \left( \theta,\lambda
\right) :0<\theta \leq \Lambda ,\,0<\lambda \leq 1\right\} \subset \mathbb{R}%
^{2}
\end{equation*}
equation \eqref{Systnew} has a unique solution in the class of
vector-functions with non-negative continuous components.
\end{lemma}

\textbf{Proof}. The proof proceeds along the lines of arguments
used to demonstrate Lemma~\ref{L_uniq}.
Starting with $\mathbf{H}^{(0)}(\theta ,\lambda):=\mathbf{C}%
_{\beta } \left( 1,\lambda \theta ^{1-\beta }\right) ^{\dagger }$ we define
\begin{eqnarray*}
\mathbf{H}^{(n+1)}(\theta ,\lambda )&:=&\mathbf{C}_{\beta }\left( 1,\lambda
\theta ^{1-\beta }\right) ^{\dagger }-\Gamma _{\beta }\theta ^{2\beta
-1}\int_0^1 \frac{\mathbf{DN}\left( \mathbf{H}^{(n)}(\theta_y,\lambda
)\right) }{\left( 1-y\right) ^{2-2\beta }}\mathrm{d}y^{\beta } \\
&=& \left( H_{1}^{(n+1)}(\theta ,\lambda ),H_{2}^{(n+1)}(\theta ,\lambda
)\right) ^{\dagger }
\end{eqnarray*}
for $n=0,1,\ldots$, where $\theta_y:=\theta \left( 1-y\right)$, $0\leq y\leq
1$.

The set $\mathcal{C}_+[0,\Lambda]\times [0,1]$ of functions $h:
[0,\Lambda]\times [0,1] \mapsto [0,\infty)\times[0,\infty)$ with continuous
components equipped with the metric
\begin{equation*}
\rho(h_1,h_2)=\sup_{\theta \in [0,\Lambda], \lambda\in [0,1]}\left\Vert
h_1(\theta, \lambda )-h_2(\theta, \lambda )\right\Vert
\end{equation*}
is a complete metric space. Invoking the Contraction Principle we conclude
that it suffices to prove
\begin{equation}  \label{inter1818}
\mathbf{H}^{(n)}\left(\theta, \lambda \right)\in \mathcal{C}%
_+[0,\Lambda]\times [0,1]
\end{equation}
for $n=0,1,\ldots$ and
\begin{equation}  \label{inter1919}
\rho(\mathbf{H}^{(n+1)}, \mathbf{H}^{(n)})\leq \kappa \rho(\mathbf{H}^{(n)},
\mathbf{H}^{(n-1)})
\end{equation}
for $n=1,2,\ldots$ and an appropriate $\kappa\in (0,1)$.

Note that all the elements of the matrices $\mathbf{D}$ and $\mathbf{C}%
_{\beta }$ are positive. Therefore, in view of the inequality
\begin{equation*}
\int_{0}^{1}\frac{\mathbf{DN}\left( \mathbf{H}^{(0)}(\theta _{y},\lambda
)\right) }{\left( 1-y\right) ^{2-2\beta }}\mathrm{d}y^{\beta }\leq
C\int_{0}^{1}\frac{\left( \lambda ^{2}\theta _{y}^{2(1-\beta )}+\lambda
\theta _{y}^{1-\beta }+1\right) }{\left( 1-y\right) ^{2-2\beta }}\mathrm{d}
y^{\beta }\mathbf{D1}\leq C\mathbf{D1}
\end{equation*}
which holds for $\beta\in (1/2,1]$, $\theta\in [0,1]$ and $\lambda\in [0,1]$%
, there exists a sufficiently small $\Lambda >0$ such that
\begin{equation*}
\mathbf{H}^{(1)}(\theta ,\lambda )=\mathbf{C}_{\beta }\left( 1,\lambda
\theta ^{1-\beta }\right) ^{\dagger }-\Gamma _{\beta }\theta ^{2\beta
-1}\int_{0}^{1}\frac{\mathbf{DN}\left( \mathbf{H}^{(0)}(\theta _{y},\lambda
)\right) }{\left( 1-y\right) ^{2-2\beta }}\mathrm{d}y^{\beta }\geq \mathbf{0}
\end{equation*}
for all $\theta \in [0,\Lambda]$ and $\lambda\in [0,1]$.

Assume now that $\mathbf{H}^{(n)}(\theta ,\lambda )\geq \mathbf{0}$ for $%
\theta\in [0,\Lambda]$ and $\lambda\in [0,1]$. Using the estimate
\begin{equation}  \label{1918}
\mathbf{H}^{(n)}(\theta ,\lambda )\leq \mathbf{C}_{\beta }\left( 1,\lambda
\theta ^{1-\beta }\right) ^{\dagger }, \ \ \theta\in [0,\Lambda], \
\lambda\in [0,1]
\end{equation}
and arguing as above we conclude $\mathbf{H}^{(n+1)}(\theta ,\lambda )\geq
\mathbf{0}$ for all $\theta\in [0,\Lambda]$ and $\lambda\in [0,1]$.
Continuity of components of $\mathbf{H}^{(n)}\left( \mathbf{\theta, \lambda }
\right)$, $n=0,1,\ldots$ is obvious, and \eqref{inter1818} follows.

Further we have
\begin{eqnarray*}
&&\left\Vert \mathbf{H}^{(n+1)}(\theta ,\lambda )-\mathbf{H}^{(n)}(\theta
,\lambda )\right\Vert \\
&&\quad \leq C_{1}\theta ^{2\beta -1}\int_{0}^{1}\left\Vert \mathbf{N}\left(
\mathbf{H}^{(n)}(\theta _{y},\lambda )\right) -\mathbf{N}\left( \mathbf{H}
^{(n-1)}(\theta _{y},\lambda )\right) \right\Vert \frac{\mathrm{d}y^{\beta }
}{\left( 1-y\right) ^{2-2\beta }} \\
&&\quad \leq C_{1}\theta ^{2\beta -1}\int_{0}^{1}\left\Vert \mathbf{H}
^{(n)}(\theta _{y},\lambda )-\mathbf{H}^{(n-1)}(\theta _{y},\lambda
)\right\Vert \frac{\mathrm{d}y^{\beta }}{\left( 1-y\right) ^{2-2\beta }} \\
&&\quad \leq C_{2}\theta ^{2\beta -1}\sup_{(x_{1},x_{2})\in \mathcal{K}%
}\left\Vert \mathbf{H}^{(n)}(x_{1},x_{2})-\mathbf{H} ^{(n-1)}(x_{1},x_{2})%
\right\Vert
\end{eqnarray*}
for all $\theta\in [0,\Lambda]$ and $\lambda\in [0,1]$, having utilized %
\eqref{NorN} and \eqref{1918} for the second inequality, and thereupon %
\eqref{inter1919} with $\kappa=C_2\Lambda ^{2\beta-1}$ for $\Lambda>0$ such
that $\kappa<1$. The proof of Lemma \ref{L_fullShur} is complete.

\noindent \textbf{Proof of Theorem \ref{T_Betamore}}. Although the
proof follows the pattern of arguments used to demonstrate Lemma
\ref{L_z3}, technical details here are more involved.

Set
\begin{equation*}
\mathbf{K}(t;\theta ,\lambda ):=\mu _{2}(t)\frac{R(t/\theta )}{\mu
_{2}(t/\theta )}\mathbf{Q}\left( t;1-\frac{\mu _{2}(t/\theta )}{\mu
_{1}R(t/\theta )},1-\frac{\lambda }{R(t/\theta )}\right)
\end{equation*}%
for positive $\theta $ and $\lambda $. It suffices to prove the
existence of a vector-function $\mathbf{K}(\theta ,\lambda )$
which, in a domain  $\mathcal{K}_{1}\subset \mathcal{K}(\beta
,\Lambda ,1)$ to be specified later, satisfies \eqref{Systnew} and
such that
\begin{equation}
\lim_{t\rightarrow \infty }\mathbf{K}(t;\theta ,\lambda )=\mathbf{K}(\theta
,\lambda ),\quad (\theta ,\lambda )\in \mathcal{K}_{1}.  \label{inter2100}
\end{equation}%
Indeed, in view of Lemma \ref{L_fullShur} $\mathbf{K}(\theta ,\lambda )=%
\mathbf{H}(\theta ,\lambda )$ for $(\theta ,\lambda )\in
\mathcal{K}_{1}$. Furthermore, using \eqref{inter2100} with
$\theta =\lambda _{1}^{1/\left( 2\beta -1\right) }$ and $\lambda
=\lambda _{2}\lambda _{1}^{-\beta /\left( 2\beta -1\right) }$ we
infer
\begin{eqnarray*}
&&\lim_{t\rightarrow \infty }R(t)\mathbf{Q}\left( t;\exp \left\{ -\lambda
_{1}\frac{\mu _{2}(t)}{\mu _{1}R(t)}\right\} ,\exp \left\{ -\frac{\lambda
_{2}}{R(t)}\right\} \right)  \\
&&\qquad=\lambda _{1}\mathbf{H}\left( \lambda _{1}^{1/\left(
2\beta -1\right) },\lambda _{2}\lambda _{1}^{-\beta /\left( 2\beta
-1\right) }\right) ,\qquad (\theta ,\lambda )\in \mathcal{K}_{1}
\end{eqnarray*}
by an appeal to Lemma \ref{L_change2} which applies because
the functions $\mu _{2}(t)/R(t)$ and $1/R(t)$ are regularly varying at $%
\infty $ with indices $1-2\beta <0$ and $-\beta <0$, respectively.

Hence we infer for $\Big(\lambda _{1}^{1/\left( 2\beta -1\right)
},\lambda _{2}\lambda _{1}^{-\beta /\left( 2\beta
-1\right)}\Big)\in \mathcal{K}_{1}$ and $i=1,2$
\begin{eqnarray}
&&\lim_{N,t\rightarrow \infty }F_{i}^{R(t)}\left( t;\exp \left\{ -\lambda
_{1}\frac{\mu _{2}(t)}{\mu _{1}R(t)}\right\} ,\exp \left\{ -\frac{\lambda
_{2}}{R(t)}\right\} \right)   \notag \\
&&\qquad\qquad=\exp \left\{ -\lambda _{1}H_{i}\left( \lambda
_{1}^{1/\left( 2\beta -1\right) },\lambda _{2}\lambda _{1}^{-\beta
/\left( 2\beta -1\right) }\right) \right\} .  \label{LALA2}
\end{eqnarray}
Since the expressions under the limits in \eqref{LALA2} are the
Laplace transforms of two-dimensional random vectors with
nonnegative components and the limits exist for
$\Big(\lambda_{1}^{1/\left( 2\beta -1\right) },\lambda _{2}\lambda
_{1}^{-\beta /\left( 2\beta -1\right)}\Big)\in \mathcal{K}_{1}$,
where $\mathcal{K}_1$ contains a ball from $(0,\infty)\times
(0,\infty)$, it follows by the uniqueness theorem for Laplace
transforms in two variables that the limit in \eqref{LALA2} exists
for all $\lambda_1>0$, $\lambda_2>0$. Moreover, there exists a
vector-function
\begin{equation*}
\mathbf{\Omega }^{\ast }\left( \lambda _{1},\lambda _{2}\right)
=\left( \Omega _{1}^{\ast }\left( \lambda _{1},\lambda _{2}\right)
,\Omega _{2}^{\ast}\left(\lambda_{1},\lambda_{2}\right)\right),\
\lambda_{1}>0,\ \lambda _{2}>0,
\end{equation*}
whose components are analytic in the domain $\left\{ \Re \lambda
_{1}>0,\Re\lambda _{2}>0\right\} $ (under an appropriate choice of
the branches) and such  that, for $i=1,2$,
\begin{equation*}
\lim_{N,t\rightarrow \infty }F_{i}^{R(t)}\left( t;\exp \left\{
-\lambda _{1}\frac{\mu _{2}(t)}{\mu _{1}R(t)}\right\} ,\exp
\left\{ -\frac{\lambda _{2}}{R(t)}\right\} \right) =e^{-\Omega
_{i}^{\ast }\left( \lambda _{1},\lambda _{2}\right)}
\end{equation*}

In particular,
\begin{equation*}
\lambda _{1}H_{i}\left( \lambda _{1}^{1/\left( 2\beta -1\right)
},\lambda _{2}\lambda _{1}^{-\beta /\left( 2\beta -1\right)
}\right) =\Omega _{i}^{\ast }\left( \lambda _{1},\lambda
_{2}\right)
\end{equation*}
for $\Big(\lambda _{1}^{1/\left( 2\beta -1\right) },\lambda
_{2}\lambda _{1}^{-\beta /\left( 2\beta -1\right) }\Big)\in
\mathcal{K}_{1}$. Since $\lambda_1^{-1}\mathbf{\Omega }^{\ast
}\left( \lambda _{1},\lambda _{2}\right) $ solves (\ref{Systnew})
for $\Big(\lambda _{1}^{1/\left( 2\beta -1\right) },\lambda
_{2}\lambda _{1}^{-\beta /\left( 2\beta -1\right) }\Big)\in
\mathcal{K}_{1}$ it follows from the uniqueness theorem for
analytic functions in two variables that
$\lambda_1^{-1}\mathbf{\Omega }^{\ast }\left( \lambda _{1},\lambda
_{2}\right) $ solves \eqref{Systnew} in the domain $\left\{
\Re\lambda _{1}>0,\Re\lambda _{2}>0\right\} $ and, therefore,
coincides with $\mathbf{H}\left( \lambda _{1}^{1/\left( 2\beta
-1\right) },\lambda _{2}\lambda _{1}^{-\beta /\left( 2\beta
-1\right) }\right) $ as desired.

The major step towards proving \eqref{inter2100} is to check that
\begin{equation}
\lim_{t\rightarrow \infty }\sup_{(\theta, \lambda) \in
\mathcal{K}_{1}}\left\Vert \mathbf{K}(t;\theta ,\lambda
)-\mathbf{H}\left( \theta ,\lambda \right) \right\Vert =0.
\label{inter2222}
\end{equation}%
To this end, we use \eqref{RenewalQ} and the equality
\begin{equation*}
\mathbf{Q}\left( t_{y};1-\frac{\mu _{2}(t/\theta )}{\mu _{1}R(t/\theta )},1-%
\frac{\lambda }{R(t/\theta )}\right) =\frac{\mu _{2}(t/\theta )}{R(t/\theta )%
}\frac{1}{\mu _{2}(t_{y})}\mathbf{K}(t_{y};\theta _{y},\lambda ),\ \ y\in
\lbrack 0,1)
\end{equation*}%
where $t_{y}=t(1-y)$ and $\theta _{y}=\theta (1-y)$, to obtain
\begin{eqnarray}
\mathbf{K}(t;\theta ,\lambda ) &=&{\frac{\mu _{2}(t)}{\mu _{2}(t/\theta )}}%
R(t/\lambda )\bigg(\mathbf{U}\ast \left( \left( \mathbf{1}-\mathbf{s}\right)
\otimes \left( \mathbf{1}-\mathbf{G}(\cdot )\right) \right) (t)  \notag
\label{Approx222} \\
&&-\int_{0}^{1}\mathrm{d}\mathbf{U}_{\mathbf{I}}(ty)\mathbf{\Phi }(\mathbf{Q}%
(t_{y};\mathbf{s}))\bigg)  \notag \\
&=&{\frac{\mu _{2}(t)}{\mu _{2}(t/\theta )}}R(t/\lambda )\bigg(\mathbf{U}%
\ast \left( \left( \mathbf{1}-\mathbf{s}\right) \otimes \left( \mathbf{1}-%
\mathbf{G}(\cdot )\right) \right) (t) \\
&&-\int_{0}^{1-\varepsilon }\mathrm{d}\mathbf{U}_{\mathbf{I}}(ty)\mathbf{%
\Phi }\bigg(\frac{\mu _{2}(t/\theta )}{R(t/\theta )}\frac{1}{\mu _{2}(t_{y})}%
\mathbf{K}(t_{y};\theta _{y},\lambda )\bigg)\bigg)-\mathbf{J}%
_{1}(t;\varepsilon ,\theta )  \notag
\end{eqnarray}%
for positive $\theta $ and $\lambda $ and any fixed $\varepsilon \in (0,1)$,
where
\begin{equation*}
\mathbf{s}=(s_{1},s_{2}):=\bigg(1-{\frac{\mu _{2}(t/\theta )}{\mu
_{1}R(t/\theta )}},1-{\frac{\lambda }{R(t/\theta )}}\bigg)
\end{equation*}%
and
\begin{equation*}
\mathbf{J}_{1}(t;\varepsilon ,\theta ):=\frac{\mu _{2}(t)}{\mu _{2}(t/\theta
)}R(t/\theta )\int_{t(1-\varepsilon )}^{t}\mathrm{d}\mathbf{U}_{\mathbf{I}%
}(w)\mathbf{\Phi }(\mathbf{Q}(t-w;\mathbf{s})).
\end{equation*}

Noting that $\mu_2(t)$ is regularly varying at $\infty$ with index $1-\beta$
and using \eqref{AsMat} and then Lemma \ref{L_firsrMom} give
\begin{eqnarray*}
\frac{\mu _{2}(t)}{\mu _{2}(t/\theta )}R(t/\theta )\mathbf{U}\ast \left(
\left( \mathbf{1-s}\right) \otimes \left( \mathbf{1-G}(\cdot \right) \right)
(t) =\frac{\mu _{2}(t)}{\mu _{2}(t/\theta )}R(t/\theta )\mathbf{P}(t)(
\mathbf{1-s}) \\
=\left( 1+\delta _{0}(t;\theta ,\lambda )\right) \mathbf{C}_{\beta }\left(
1,\lambda \theta ^{1-\beta }\right) ^{\dagger }
\end{eqnarray*}
where $\lim_{t\to\infty}\delta_0(t;\delta, \lambda)=0$ uniformly in $%
(\delta, \lambda) \in\mathcal{K}$ (while uniformity in $\theta$ follows from
the regular variation of $\mu_2(t)$ and Theorem 1.5.2 in \cite{BGT},
uniformity in $\lambda$ is trivial).

To proceed we intend to show that
\begin{equation}  \label{inter24}
\int_{t(1-\varepsilon )}^t{\frac{\mathrm{d}\mathbf{U}_{\mathbf{I}}(y)}{\hat{%
\mu}_2^2(t-y)}}\mathbf{1} \ \sim \ {\frac{R(t)}{\mu^2_2(t)}}%
\int_{1-\varepsilon}^1(1-y)^{2\beta-2}y^{\beta-1}\mathrm{d}y\beta\Gamma_\beta%
\mathbf{D}\mathbf{1}, \ \ t\to\infty
\end{equation}
for any fixed $\varepsilon\in (0,1)$, where, as before, $\hat{\mu}%
_{2}(t)=\mu _{2}(t)+1$.

Since
\begin{equation}  \label{lim_U1}
\lim_{t\to\infty}\mathbf{U}_{\mathbf{I}}(ty)/R(t)=\Gamma_\beta\mathbf{D}%
y^\beta, \ \ y>0
\end{equation}
by Lemma \ref{L_asympRen} and $\lim_{t\to\infty}\hat{\mu}^2_{2}(t)/ \hat{\mu}%
^2_{2}(t(1-y))=(1-y)^{2\beta-2}$ uniformly in $y\in [0,1-\varepsilon]$ we
conclude that
\begin{equation}  \label{inter25}
\lim_{t\to\infty}\int_0^{1-\varepsilon}{\frac{\hat{\mu}^2_{2}(t)}{\hat{\mu}%
^2_2(t(1-y))}}{\frac{\mathrm{d}\mathbf{U}_{\mathbf{I}}(ty)}{R(t)}}=
\int_0^{1-\varepsilon}(1-y)^{2\beta-2}y^{\beta-1}\mathrm{d}y\beta\Gamma_\beta%
\mathbf{D}
\end{equation}
by an application of Lemma A.5 in \cite{Iksanov:2013}. A combination of %
\eqref{inter25} with
\begin{equation*}
\int_0^1{\frac{\mathrm{d}\mathbf{U}_{\mathbf{I}}(y)}{\hat{\mu}_2^2(t-y)}} \
\sim \ {\frac{R(t)}{\mu^2_2(t)}}\int_0^1(1-y)^{2\beta-2}y^{\beta-1}\mathrm{d}%
y\beta\Gamma_\beta\mathbf{D}, \ \ t\to\infty
\end{equation*}
which is just an equivalent form of one part of Lemma \ref{CShur2} proves %
\eqref{inter24}.

As another preparatory result we obtain
\begin{eqnarray}
\mathbf{Q}(w,\mathbf{s}) &\leq &\mathbf{P}(w)\left( \mathbf{1}-\mathbf{s}
\right) =\frac{1}{R(t/\theta )}\mathbf{P}(w)\left(
\begin{array}{c}
\mu _{2}(t/\theta )\mu _{1}^{-1} \\
\lambda%
\end{array}
\right)  \notag \\
&\leq &\frac{C}{R(t/\theta )}\left( \frac{\mu _{2}(t/\theta )}{\hat{\mu}
_{2}(w)}+\lambda \right) \mathbf{1}\leq C\frac{\mu _{2}(t/\theta )}{
R(t/\theta )\hat{\mu}_{2}(w)}\mathbf{1}  \label{Star11}
\end{eqnarray}
for $w\leq t$ and positive $\theta$ and $\lambda$, where the second
inequality follows from
\begin{equation*}
\mathbf{P}_1(t)\leq C/\hat{\mu}_{2}(t), \ \ \mathbf{P}_2(t)\leq C, \ \ t\geq
0
\end{equation*}
which, in its turn, is a consequence of Lemma \ref{L_firsrMom}.

Using \eqref{Star11} and then \eqref{inter24} we have
\begin{eqnarray*}
\mathbf{J}_{1}(t;\varepsilon ,\theta ) &\leq& \frac{\mu _{2}(t)}{\mu
_{2}(t/\theta )}R(t/\theta )\int\limits_{t(1-\varepsilon )}^{t}\mathrm{d}
\mathbf{U}_{\mathbf{I}}(w)\mathbf{\Phi }\left( C\frac{\mu _{2}(t/\theta )}{%
R(t/\theta )\hat{\mu}_{2}(t-w)}\mathbf{1}\right)  \notag \\
&\leq&C\mu_2(t)\frac{\mu_{2}(t/\theta)}{R(t/\theta )}\int_{t(1-\varepsilon
)}^{t}\frac{\mathrm{d}\mathbf{U}_{\mathbf{I}}(w)}{\hat{ \mu}_{2}^{2}(t-w)}%
\mathbf{1}  \notag \\
&\sim& {\frac{\mu_2(t/\theta)}{R(t/\theta)}}{\frac{R(t)}{\mu_2(t)}}%
\int_{1-\varepsilon}^1(1-y)^{2\beta-2}y^{\beta-1}\mathrm{d}%
yC\beta\Gamma_\beta\mathbf{D}\mathbf{1}
\end{eqnarray*}
for positive $\theta$ and $\lambda$, and any fixed $\varepsilon\in (0,1)$
which implies
\begin{equation*}
\underset{\varepsilon\downarrow 0}{\lim}\,\underset{t\to\infty}{\lim\sup}
\,\sup_{(\theta, \lambda)\in\mathcal{K}}||\mathbf{J}_1(t;\varepsilon,%
\theta)||=0
\end{equation*}
by an application of Theorem 1.5.2 in \cite{BGT}.

We divide the subsequent proof into two parts according to whether $\mathbf{%
\Phi }\left( \mathbf{x}\right) =\mathbf{N}(\mathbf{x})$ or $\mathbf{\Phi}$
is not restricted in this way.

\noindent\textsc{Case $\mathbf{\Phi }=\mathbf{N}$}, i.e., the generating
functions of the reproduction laws are polynomials of degree $2$. The
advantage of this case is that equality \eqref{Approx222} takes a simpler
form
\begin{eqnarray}
&&\hspace{-0.9cm}\mathbf{K}(t;\theta ,\lambda ) =\frac{\mu _{2}(t)}{\mu
_{2}(t/\theta )}R(t/\theta )\mathbf{U}\ast \left( (\mathbf{1-s})\otimes
\left( \mathbf{1-G}(\cdot \right) \right) (t)-\mathbf{J}_{1}(t;\varepsilon
,\theta )  \notag \\
&&\qquad \quad -\frac{\mu _{2}(t/\theta )R(t)}{\mu _{2}(t)R(t/\theta )}
\int_{0}^{1-\varepsilon }\frac{\mathrm{d}\mathbf{U}_{\mathbf{I}}(ty)}{R(t)}
\frac{\mu _{2}^{2}(t)}{\mu _{2}^{2}(t_{y})}\mathbf{N}\left( \mathbf{K}
(t_{y};\theta _{y},\lambda )\right)  \label{Approx2222}
\end{eqnarray}
which, in view of our findings in the preceding part of the proof, can be
represented as follows
\begin{eqnarray*}
\mathbf{K}(t;\theta ,\lambda ) &=&\mathbf{C}_{\beta }\left( 1,\lambda \theta
^{1-\beta }\right) ^{\dagger }+\delta _{1}(t;\varepsilon ,\theta ,\lambda )
\mathbf{1} \\
&&-\theta ^{2\beta -1}\int_{0}^{1-\varepsilon }\frac{\mathrm{d}\mathbf{U}_{
\mathbf{I} }(ty)}{R(t)}\frac{\mathbf{N}\left( \mathbf{K}(t_{y};\theta
_{y},\lambda )\right) }{\left( 1-y\right) ^{2-2\beta }}
\end{eqnarray*}
where
\begin{equation}
\lim_{\varepsilon \downarrow 0}\mathop{\overline{\lim}}\limits_{t\rightarrow
\infty }\sup_{\left( \theta ,\lambda \right) \in \mathcal{K}}\left\vert
\delta _{1}\left( t;\varepsilon ,\theta ,\lambda \right) \right\vert =0.
\label{Del}
\end{equation}
Applying \eqref{lim_U1} and arguing as in the proof of Lemma \ref{L_z3} we
can show that
\begin{eqnarray*}
\mathbf{H}(\theta ,\lambda ) &=&\mathbf{C}_{\beta }\left( 1,\lambda \theta
^{1-\beta }\right) ^{\dagger }+\delta _{2}\left( t;\varepsilon ,\theta
,\lambda \right) \mathbf{1} \\
&&-\theta ^{2\beta -1}\int_{0}^{1-\varepsilon }\frac{\mathrm{d}\mathbf{U}_{
\mathbf{I} }(ty)}{R(t)}\frac{\mathbf{N}\left( \mathbf{H}(\theta _{y},\lambda
)\right) }{\left( 1-y\right) ^{2-2\beta }}
\end{eqnarray*}
for $(\theta,\lambda)\in\mathcal{K}$ and any fixed $\varepsilon\in (0,1)$,
where
\begin{equation}
\lim_{\varepsilon \downarrow 0}\mathop{\overline{\lim}}\limits_{t\rightarrow
\infty }\sup_{\left( \theta ,\lambda \right) \in \mathcal{K}}\left\vert
\delta _{2}\left( t;\varepsilon ,\theta ,\lambda \right) \right\vert =0.
\label{Rho}
\end{equation}

Put
\begin{eqnarray*}
\mathbf{V}\left( t_{y};\theta ,\lambda \right) &:= &\mathbf{N}(\mathbf{K}
(t_{y};\theta ,\lambda ))-\mathbf{N}\left( \mathbf{H}(\theta ,\lambda
)\right) , \\
&& \\
\mathbf{L}\left( t;\theta ,\lambda \right) &:= &\left\Vert \mathbf{K}
(t;\theta ,\lambda )-\mathbf{H}(\theta ,\lambda )\right\Vert .
\end{eqnarray*}
Using \eqref{Star11} with $w=t_y$ for $y\in(0,1-\varepsilon)$ gives
\begin{equation*}
\mathbf{K}(t_{y};\theta _{y},\lambda )\leq C\mathbf{1}
\end{equation*}
which together with the boundedness of $\mathbf{H}(\theta ,\lambda )$ in $%
\left( \theta ,\lambda \right) \in \mathcal{K}$ which follows from %
\eqref{Systnew} and \eqref{NorN} implies that
\begin{equation*}
\left\Vert \mathbf{V}\left( t_{y};\theta _{y},\lambda \right) \right\Vert
\leq C\mathbf{L}(t_{y};\theta _{y},\lambda ).
\end{equation*}
for $\left( \theta ,\lambda \right) \in \mathcal{K}$. Setting $\delta
(t;\varepsilon ,\theta ,\lambda ):=\left\vert \delta _{1}(t;\varepsilon
,\theta ,\lambda )\right\vert +\left\vert \delta _{2}(t;\varepsilon ,\theta
,\lambda )\right\vert $ we get
\begin{eqnarray*}
\mathbf{L}\left( t;\theta ,\lambda \right) &\leq &C\delta (t;\varepsilon
,\theta ,\lambda ) \\
&&\mathbf{+}C\theta ^{2\beta -1}\left\Vert \int_{0}^{1-\varepsilon }\frac{
\mathrm{d} \mathbf{U}_{\mathbf{I}}(tw)}{R(t)}\frac{\mathbf{V}\left(
t_{y};\theta _{y},\lambda \right) }{\left( 1-y\right) ^{2-2\beta }}
\right\Vert \\
&\leq &C\delta (t;\varepsilon ,\theta ,\lambda )\mathbf{+}C\theta ^{2\beta
-1}\int_{0}^{1-\varepsilon }\left\Vert \frac{\mathrm{d}\mathbf{U}_{\mathbf{I}
}(tw)}{ R(t)}\right\Vert \frac{\mathbf{L}(t_{y};\theta _{y},\lambda )}{
\left( 1-y\right) ^{2-2\beta }}
\end{eqnarray*}
and thereupon
\begin{equation*}
\mathbf{L}\left( t;\theta ,\lambda \right) \leq C\delta (t;\varepsilon
,\theta ,\lambda )+C_{1}\theta ^{2\beta -1}\sup_{w\geq t\varepsilon
}\sup_{\left( \theta ^{\ast },\lambda ^{\ast }\right) \in \mathcal{K}}
\mathbf{L}(w;\theta ^{\ast },\lambda ^{\ast })  \label{Imp1}
\end{equation*}
for $\left( \theta ,\lambda \right) \in \mathcal{K}$. By shrinking, if
needed, $\mathcal{K}$ (which gives $\mathcal{K}_1$) we can and do assume
that $\kappa:=C_1\Lambda^{2\beta-1}<1$. Setting
\begin{equation*}
S(t):=\sup_{w\geq t}\sup_{\left( \theta ,\lambda \right) \in \mathcal{K}_1}
\mathbf{L}(w;\theta ,\lambda )
\end{equation*}
and letting first $t\to\infty$ and then $\varepsilon\downarrow 0$ (at which
step relations \eqref{Del} and \eqref{Rho} have to be recalled) in the
inequality
\begin{equation*}
S(t)\leq C\sup_{s\geq t}\sup_{(\theta,\lambda)\in\mathcal{K}_1}\delta
(t;\varepsilon ,\theta ,\lambda )+\kappa S(t\varepsilon )
\end{equation*}
we arrive at \eqref{inter2222}.

This completes the proof of the theorem under the assumption $\mathbf{\Phi }=%
\mathbf{N}$.

\noindent \textsc{General case} can be treated along the same lines after
noting that representation \eqref{DDEF} along with relations \eqref{Star11}
and \eqref{inter25} implies
\begin{eqnarray*}
&&\frac{\mu _{2}(t)}{\mu _{2}(t/\theta )}R(t/\theta )\left\Vert
\int_{0}^{1-\varepsilon }\mathrm{d}\mathbf{U}_{\mathbf{I}}(ty)\left[ \mathbf{%
\ \Phi }(\mathbf{Q}(t_{y};\mathbf{s}))-\mathbf{N}\left( \mathbf{Q}(t_{y};
\mathbf{s})\right) \right] \right\Vert \\
&&\qquad \qquad \qquad =o\left( \frac{\mu _{2}(t)}{\mu _{2}(t/\theta )}
R(t/\theta )\int_{0}^{1-\varepsilon }\left\Vert \mathrm{d}\mathbf{U}_{
\mathbf{I}}(ty)\right\Vert \left\Vert \frac{\mu _{2}(t/\theta )}{R(t/\theta
) \hat{\mu}_{2}(t_{y})}\mathbf{1}\right\Vert ^{2}\right) \\
&&\qquad \qquad \qquad =o\left( {\frac{\mu_2(t)\mu _{2}(t/\theta )}{%
R(t/\theta )}}\int_0^{1- \varepsilon }\frac{\left\Vert \mathrm{d}\mathbf{U}_{%
\mathbf{I}}(ty)\right\Vert }{\hat{\mu}_2^2(t_y)}\right) \\
&&\qquad \qquad \qquad =o\left( {\frac{\mu_2(t/\theta)}{R(t/\theta)}}{\frac{%
R(t)}{\mu_2(t)}} \right) =o\left( 1\right)
\end{eqnarray*}
for any $\varepsilon \in\left( 0,1\right)$. Furthermore, this convergence is
uniform in $\left( \theta ,\lambda \right) \in \mathcal{K}_1$ by Theorem
1.5.2 in \cite{BGT}. The proof of Theorem \ref{T_Betamore} is complete.

\section{Proofs for the final evolutionary stages\label{Sec5}}

The purpose of this section is to prove Theorem \ref{T_zone2} which will be
done in a series of lemmas.

The following result was obtained in Lemma 2 in \cite{13V}.

\begin{lemma}
\label{L_growth} Assume that the Perron root of the aperiodic irreducible
matrix $\mathbf{M}$ in $(\ref{01})$ is equal to $1$. Then the generating
vector-function $\mathbf{F}(t;\mathbf{s})$ is non-decreasing in $t$ for each
$\mathbf{s}\in \mathcal{A}:=\left\{ \mathbf{s}\in \left[ 0,1 \right] ^{2}:%
\mathbf{f}\left( \mathbf{s}\right) \geq \mathbf{s}\right\}$.
\end{lemma}

Recalling (\ref{DefMatrix}) we set
\begin{equation*}
{\mathfrak{Q}}(t;\mathbf{s}):=(\mathbf{v},\mathbf{Q}(t;\mathbf{s})).
\end{equation*}
The next statement is a natural generalization of Theorem 1 in \cite{Chi71}.

\begin{lemma}
\label{L_ratio} Let $\mathbf{Z}(t)$ be an irreducible, aperiodic, and
critical process satisfying conditions $(\ref{02})$ and $(\ref{03})$. Then,
for each $\lambda >0$ and any non-increasing function $\psi$ such that $%
\lim_{t\to\infty}\psi (t)=0$, we have
\begin{equation*}
\lim_{t\rightarrow \infty }\frac{\mathbf{1}-\mathbf{F}(t;\mathbf{1}-\lambda
\mathbf{u}\psi (t))}{{\mathfrak{Q}}(t;\mathbf{1}-\lambda \mathbf{u}\psi (t))}
=\lim_{t\rightarrow \infty }\frac{\mathbf{Q}(t;\mathbf{1}-\lambda \mathbf{u}
\psi (t))}{{\mathfrak{Q}}(t;\mathbf{1}-\lambda \mathbf{u}\psi (t))}=\mathbf{u%
}.
\end{equation*}
\end{lemma}

\noindent \textbf{Proof}. For $t$ large enough to ensure $\mathbf{1}-\lambda
\mathbf{u}\psi (t)\in \left[ 0,1\right]^2$ use
\begin{equation*}
\mathbf{1}-\mathbf{f}\left( \mathbf{1}-\lambda \mathbf{u}\psi (t)\right)
\leq \lambda \psi (t)\mathbf{Mu=}\lambda \psi (t)\mathbf{u}
\end{equation*}
to infer $\mathbf{1}-\lambda \mathbf{u}\psi (t)\in \mathcal{A}$ and further
\begin{eqnarray*}
&&\mathbf{F}(t;\mathbf{1}-\lambda \mathbf{u}\psi (t))\geq \mathbf{F}(t-w;
\mathbf{1}-\lambda \mathbf{u}\psi (t)) \\
&\geq& \mathbf{F}(t-w;\max (1-\lambda u_1\psi(t-w),0), \max(1-\lambda
u_2\psi(t-w),0))
\end{eqnarray*}
for each $w\in (0,t]$, having utilized Lemma \ref{L_growth} for the first
inequality and the monotonicity of $\psi$ for the second. A minor
modification of the proof of Theorem 1 in \cite{Chi71} finishes the proof of
the lemma.

\begin{lemma}
\label{L_Qsmall} Suppose that Hypothesis \textbf{A} holds. Then, for each $%
\lambda >0$ and any non-increasing function $\psi$ such that $%
\lim_{t\to\infty}\psi (t)=0$ and $1/R(t)=o(\psi (t))$, we have
\begin{equation*}
{\mathfrak{Q}}(t;\mathbf{1}-\lambda \mathbf{u}\psi (t))=o(\psi (t)).
\end{equation*}
\end{lemma}

\noindent \textbf{Proof}. In view of \eqref{Lu} and Lemma \ref{L_asympRen}
\begin{equation}
\mathbf{U}_{\mathbf{I}}(t)\sim \mathbf{U}(t)\sim \Gamma _{\beta }R(t)
\mathbf{\ D}.  \label{DefU1}
\end{equation}
According to Lemma \ref{L_growth} the vector-function $\mathbf{\Phi }(
\mathbf{Q}(t;\mathbf{s}))$ is non-increasing in $t$ for each fixed $\mathbf{%
s }\in \mathcal{A}$. This enables us to infer from \eqref{RenewalQ}
\begin{eqnarray}
\mathbf{U}_{\mathbf{I}}(t)\mathbf{\Phi }(\mathbf{Q}(t;\mathbf{s})) &\leq
&\int_{0}^{t}\mathrm{d}\mathbf{U}_{\mathbf{I}}(w)\mathbf{\Phi }(\mathbf{Q}
(t-w;\mathbf{s}))  \notag \\
\qquad \qquad &\leq &\mathbf{U}\ast \left( \left( \mathbf{1}-\mathbf{s}
\right) \otimes \left( \mathbf{1}-\mathbf{G}(\cdot )\right) \right) (t)
\label{inter}
\end{eqnarray}
for each fixed $\mathbf{s}\in \mathcal{A}$.

Put $\mathbf{s}=\mathbf{s}(t):=\mathbf{1}-\lambda \mathbf{u}\psi (t)$ and
use \eqref{AsMat} and Lemma \ref{L_firsrMom} to obtain
\begin{equation*}
\mathbf{U}\ast \left( \left( \mathbf{1}-\mathbf{s}\right) \otimes \left(
\mathbf{1}-\mathbf{G}(\cdot )\right) \right) (t)=\left( 1+o(1)\right)
\lambda \psi (t)\mathbf{DJ}(t)\mathbf{1}=O(\psi (t))\mathbf{1}
\end{equation*}
which entails
\begin{equation*}
\mathbf{U}_{\mathbf{I}}(t)\mathbf{\Phi }(\mathbf{Q}(t;\mathbf{1}-\lambda
\mathbf{u}\psi (t)))=O(\psi (t))\mathbf{1}
\end{equation*}
by an appeal to \eqref{inter} which applies because $\mathbf{s}(t)\in
\mathcal{A}$ for large enough $t$ (see the proof of Lemma \ref{L_ratio}).
Recalling \eqref{DefU1} we conclude
\begin{equation*}
\mathbf{\Phi }(\mathbf{Q}(t;\mathbf{1}-\lambda \mathbf{u}\psi (t)))=O\left(
\psi (t)/R(t)\right) \mathbf{1}=o\left( \psi ^{2}(t)\right) \mathbf{1}
\end{equation*}
which in view of (\ref{DDEF}) and (\ref{DefN}) ensures
\begin{equation*}
\mathbf{Q}(t;\mathbf{1}-\lambda \mathbf{u}\psi (t))=o\left( \psi (t)\right)
\mathbf{1}.
\end{equation*}
The proof of Lemma \ref{L_Qsmall} is complete.

\begin{lemma}
\label{L_z2} Suppose that Hypothesis \textbf{A}$(0,1)$ and conditions %
\eqref{Nzone1} hold, and that $1/R(t)=o(\psi (t))$ for a non-increasing
regularly varying function $\psi (t)=t^{-\gamma }\ell_2(t)$, $\gamma\in
[0,1) $, such that $\lim_{t\to\infty}\ell_2(t)=0$. Then, for each $\lambda
>0 $,
\begin{equation}
\mathbf{Q}(t;\mathbf{e}^{-\lambda \mathbf{u}\psi (t)})\sim \sqrt{\lambda }%
\sqrt{\frac{v_{2}u_{2}}{B}\psi (t)(1-G_2(t)) }\mathbf{u}^{\mathbf{\dagger }}.
\label{0M2}
\end{equation}
If, in addition, condition \eqref{M1} holds, then
\begin{equation}
\lim_{N,t\rightarrow \infty }N\mathbf{Q}(t;1,e^{-\lambda u_{2}\psi (t)})=r%
\sqrt{\lambda }\mathbf{u}^{\mathbf{\dagger }}.  \label{Fpart}
\end{equation}
\end{lemma}

\noindent\textbf{Proof}. Using \eqref{BBa1} with $\mathbf{s}=\mathbf{1}
-\lambda \mathbf{u}\psi (t)$ for large enough $t$ to ensure $\mathbf{s}\in
\mathcal{A}$ we get
\begin{eqnarray*}
&&\mathbf{Q}(t;\mathbf{1}-\lambda \mathbf{u}\psi (t))=\lambda \psi (t)
\mathbf{u}\otimes \left( \mathbf{1}-\mathbf{G}(t)\right) \\
&&\qquad +\int_{0}^{t}\left( \mathbf{1}-\mathbf{f}(\mathbf{F}(t-w;\mathbf{1}
-\lambda \mathbf{u}\psi (t)))\right) \otimes \mathrm{d}\mathbf{G}(w) \\
&&\quad \geq \lambda \psi (t)\mathbf{u}\otimes \left( \mathbf{1}-\mathbf{G}
(t)\right) +\left( \mathbf{1}-\mathbf{f}(\mathbf{F}(t;\mathbf{1}-\lambda
\mathbf{u}\psi (t)))\right) \otimes \mathbf{G}(t)
\end{eqnarray*}
and, therefore,
\begin{eqnarray*}
\mathbf{\Phi }(\mathbf{Q}(t;\mathbf{1}-\lambda \mathbf{u}\psi (t)))\otimes
\mathbf{G}(t) &\geq &\lambda \psi (t)\mathbf{u}\otimes \left( \mathbf{1}-
\mathbf{G}(t)\right) \\
&&\hspace{-2.0cm} +\mathbf{MQ}(t;\mathbf{1}-\lambda \mathbf{u}\psi
(t))\otimes \mathbf{G}(t)- \mathbf{Q}(t;\mathbf{1}-\lambda \mathbf{u}\psi
(t)).
\end{eqnarray*}
Multiplying both sides of this inequality by $\mathbf{v}$, the left
eigenvector of $\mathbf{M}$ corresponding to the Perron root $1$, and using
the equality $\left( \mathbf{v},\mathbf{Q}(t;\mathbf{s})\right) =\left(
\mathbf{v},\mathbf{MQ}(t;\mathbf{s})\right) $ we obtain
\begin{eqnarray}  \label{LF}
(\mathbf{v,\Phi }(\mathbf{Q}(t;\mathbf{1}-\lambda \mathbf{u}\psi (t))))
&\geq &\lambda \psi (t)(\mathbf{v},\mathbf{u}\otimes \left( \mathbf{1}-
\mathbf{G}(t)\right) ) \\
&& -\left( \mathbf{v},\mathbf{MQ}(t;\mathbf{1}-\lambda \mathbf{u}\psi
(t))\otimes (\mathbf{1}-\mathbf{G}(t))\right)  \notag
\end{eqnarray}
which in combination with
\begin{equation*}
(\mathbf{v},\mathbf{MQ}(t;\mathbf{1}-\lambda \mathbf{u}\psi (t)))=o\left(
\psi (t)\right)
\end{equation*}
(see Lemma \ref{L_Qsmall}) implies
\begin{equation}
\mathop{\lim\inf}\limits_{t\rightarrow \infty }\frac{(\mathbf{v,\Phi }(
\mathbf{Q}(t;\mathbf{1}-\lambda \mathbf{u}\psi (t))))}{\psi (t)(\mathbf{v},
\mathbf{u}\otimes \left( \mathbf{1}-\mathbf{G}(t)\right) )}\geq \lambda .
\label{LamBelow}
\end{equation}
Hence
\begin{equation}  \label{Fbel}
\mathop{\lim\inf}\limits_{T\rightarrow \infty }\frac{\textstyle{\int_{0}^{T}
t(\mathbf{v,\Phi }(\mathbf{Q}(t;\mathbf{1}-\lambda \mathbf{u}\psi (t))))
\mathrm{d}t}}{\textstyle{\ \int_{0}^{T}{t}\psi (t)(\mathbf{v},\mathbf{u}
\otimes \left( \mathbf{1}- \mathbf{G}(t)\right) )\mathrm{d}t}}\geq \lambda
\end{equation}
by a variant of l'H\^{o}pital rule which is applicable because
\begin{equation*}
\int_0^\infty t\psi (t)(\mathbf{v},\mathbf{u}\otimes \left( \mathbf{1} -
\mathbf{G}(t)\right) )\mathrm{d}t=\infty.
\end{equation*}

Now we intend to prove the converse inequality for the upper limit. To this
end, we will use the equality
\begin{eqnarray}
&&\bigg(\mathbf{v}, \int_0^t\left( \mathbf{1}-\mathbf{f}(\mathbf{F}(t-w;
\mathbf{s}))\right) \otimes \mathrm{d}\mathbf{G}(w)\bigg)-\left(\mathbf{v},
\mathbf{1}-\mathbf{f}(\mathbf{F}(t;\mathbf{s}))\right)  \notag \\
&&\qquad \qquad \qquad =(\mathbf{v}, \mathbf{\Phi }(\mathbf{Q}(t;\mathbf{s}
)))-(\mathbf{v}, \left( \mathbf{1}-\mathbf{s}\right) \otimes \left( \mathbf{%
1 }-\mathbf{G}(t)\right))  \label{Int1}
\end{eqnarray}
which is a consequence of \eqref{BBa1} and $\mathbf{v}\mathbf{M}=\mathbf{v}$%
. Denoting
\begin{eqnarray*}
\mathbf{\Delta }(T) &:=&\int_0^T t\mathrm{d}t\int_0^t\left( \mathbf{1}-
\mathbf{f}(\mathbf{F}(t-w;\mathbf{s}(t)))\right) \otimes \mathrm{d}\mathbf{G}
(w) \\
&&-\int_0^T t\left( \mathbf{1}-\mathbf{f}(\mathbf{F}(t;\mathbf{s}
(t)))\right) \mathrm{d}t \\
&=&\int_0^T \mathrm{d}\mathbf{G}(w)\otimes \int_0^{T-w}\left( t+w\right)
\left( \mathbf{1}-\mathbf{f}(\mathbf{F}(t;\mathbf{s}(t+w)))\right) \mathrm{d}
t \\
&&-\int_{0}^{T}t\left( \mathbf{1}-\mathbf{f}(\mathbf{F}(t;\mathbf{s}
(t)))\right) \mathrm{d}t,
\end{eqnarray*}
where $\mathbf{s}(t)=\mathbf{1}-\lambda \mathbf{u}\psi (t)$ we observe that %
\eqref{Int1} with $\mathbf{s}=\mathbf{s}(t)$ is equivalent to
\begin{equation}  \label{inter7}
(\mathbf{v}, \mathbf{\Delta }(T))=\int_0^T t(\mathbf{v}, \mathbf{\Phi }(
\mathbf{Q}(t;\mathbf{s}(t))))\mathrm{d}t-\lambda \int_0^T t\psi (t)(\mathbf{%
v }, \mathbf{u}\otimes \left( \mathbf{1}-\mathbf{G}(t)\right)) \mathrm{d}t.
\end{equation}

By the monotonicity of $\psi(t)$ we have
\begin{equation*}
\int_0^T \mathrm{d}\mathbf{G}(w)\otimes \int_0^{T-w} t\left( \mathbf{1}-
\mathbf{f}(\mathbf{F}(t;\mathbf{s}(t+w)))\right) \mathrm{d}t\leq \int_0^T
t\left( \mathbf{1}-\mathbf{f}(\mathbf{F}(t;\mathbf{s}(t)))\right)\mathrm{d}t
\end{equation*}
which implies
\begin{eqnarray*}
\mathbf{\Delta }(T) &\leq &\int_0^T \mathrm{d}\mathbf{G}(w)\otimes
\int_{0}^{T-w}\left[ \left( t+w\right) -t\right] \left( \mathbf{1}-\mathbf{f}
(\mathbf{F}(t;\mathbf{s}(t+w)))\right) \mathrm{d}t \\
&=&\int_{0}^T w \mathrm{d}\mathbf{G}(w)\otimes \int_{0}^{T-w}\left( \mathbf{%
1 }-\mathbf{f}(\mathbf{F}(t;\mathbf{s}(t+w)))\right)\mathrm{d}t
\end{eqnarray*}
and further
\begin{eqnarray*}
\mathbf{\Delta }(T) &\leq &\int_0^T w\mathrm{d}\mathbf{G}(w)\otimes
\int_0^{T-w}\mathbf{MQ}(t;\mathbf{s}(t+w)))\mathrm{d}t \\
&\leq &\int_0^T w\mathrm{d}\mathbf{G}(w)\otimes \int_0^T \mathbf{MQ}(t;
\mathbf{s}(t)))\mathrm{d}t \\
&=&\int_0^T w\mathrm{d}\mathbf{G}(w)\otimes \int_{0}^{T}o\left( \psi (t)
\mathbf{1}\right) \mathrm{d}t=o\bigg(\mathbf{1} \int_0^T w\int_0^T\psi (t)
\mathrm{d}t\mathrm{d}\mathbf{G}(w)\bigg),
\end{eqnarray*}
having utilized $\mathbf{\Phi}(\mathbf{Q}(t;\mathbf{s}))=\mathbf{M}\mathbf{Q}
(t;\mathbf{s})-(\mathbf{1}-\mathbf{f}(\mathbf{F}(t;\mathbf{s})))\geq \mathbf{%
\ 0}$ for the first inequality, and the monotonicity of $\psi$ for the
second. The first equality follows from Lemmas \ref{L_ratio} and \ref%
{L_Qsmall}, and the second equality is a consequence of $\int_0^\infty \psi
(t)\mathrm{d} t=\infty$.

Recalling that $\psi$ is regularly varying at $\infty$ with index $-\gamma$,
$\gamma \in [0,1)$ and invoking Proposition 1.5.8 in \cite{BGT} we infer
\begin{equation*}
\int_0^T w\mathrm{d}G_2(w)\int_0^T\psi(t)\mathrm{d}t \sim {\frac{\beta}{
(1-\beta)(1-\gamma)}}T^2\psi(T)(1-G_2(T))
\end{equation*}
and
\begin{equation}  \label{inter8}
\int_0^T t \psi (t)(\mathbf{v},\mathbf{u}\otimes \left( \mathbf{1}-\mathbf{G}
(t)\right) )\mathrm{d}t \sim {\frac{v_2u_2}{2-\gamma-\beta}}
T^2\psi(T)(1-G_2(T))
\end{equation}
which implies
\begin{equation*}
\int_0^T w\mathrm{d}\mathbf{G}(w)\int_0^T \psi (t)\mathrm{d}t=O\left(
\mathbf{1} \int_{0}^{T}{t}\psi (t)(\mathbf{v},\mathbf{u}\otimes \left(
\mathbf{1}-\mathbf{G}(t)\right) )\mathrm{d}t\right)
\end{equation*}
for $\int_0^\infty w\mathrm{d}G_1(w)=\mu_1<\infty$. Therefore
\begin{equation*}
\left( \mathbf{v},\mathbf{\Delta }(T)\right)=o\bigg(\int_0^T t\psi (t)(
\mathbf{v},\mathbf{u}\otimes \left( \mathbf{1}-\mathbf{G}(t)\right) )\mathrm{%
\ d}t\bigg)
\end{equation*}
which entails
\begin{equation*}
\mathop{\lim\sup}\limits_{T\to \infty}\frac{\int_0^T t(\mathbf{v}, \mathbf{\
\Phi }(\mathbf{Q} (t;\mathbf{s}(t))))\mathrm{d}t}{\int_0^T t\psi (t)(\mathbf{%
\ v}, \mathbf{u}\otimes \left( \mathbf{1 }-\mathbf{G}(t)\right)) \mathrm{d}t}
\leq \lambda
\end{equation*}
in view of \eqref{inter7}. Combining this with \eqref{Fbel} gives
\begin{equation*}
\lambda \int_0^T t\psi (t)(\mathbf{v},\mathbf{u}\otimes \left( \mathbf{1}-
\mathbf{G}(t)\right) )\mathrm{d}t \sim \int_0^T t(\mathbf{v,\Phi }(\mathbf{Q}
(t;\mathbf{1}-\lambda \mathbf{u}\psi (t))))\mathrm{d}t.
\end{equation*}
Left with the proof of
\begin{equation}  \label{inter9}
\lambda \psi (t)(\mathbf{v},\mathbf{u}\otimes \left( \mathbf{1}-\mathbf{G}
(t)\right) )\sim (\mathbf{v,\Phi }(\mathbf{Q}(t;\mathbf{1}-\lambda \mathbf{u}
\psi (t))))
\end{equation}
we infer by means of \eqref{inter8}
\begin{equation*}
\int_T^{\delta T}t(\mathbf{v,\Phi }(\mathbf{Q}(t;\mathbf{1}-\lambda \mathbf{%
u }\psi (t)))) \mathrm{d}t \sim {\frac{\lambda(\delta^{2-\gamma-\beta}-1) }{
2-\gamma-\beta}}T^2\psi(T)(\mathbf{v},\mathbf{u}\otimes \left( \mathbf{1}-
\mathbf{G}(T)\right) )
\end{equation*}
for $\delta>1$. The function $t\mapsto (\mathbf{v,\Phi } (\mathbf{Q}(t;%
\mathbf{1}-\lambda \mathbf{u}\psi (t))))$ is eventually non-increasing
whence
\begin{equation*}
\int_T^{\delta T}t(\mathbf{v,\Phi }(\mathbf{Q}(t;\mathbf{1}-\lambda \mathbf{%
u }\psi (t)))) \mathrm{d}t\leq (\delta-1)\delta T^2(\mathbf{v,\Phi }(\mathbf{%
Q} (T;\mathbf{1}-\lambda \mathbf{u}\psi(T))))
\end{equation*}
for large enough $T$, and thereupon
\begin{equation*}
\underset{t\to\infty}{\lim\inf}\,{\frac{(\mathbf{v,\Phi }(\mathbf{Q}(t;
\mathbf{1}-\lambda \mathbf{u}\psi (t))))}{\psi(t)(\mathbf{v},\mathbf{u}
\otimes \left( \mathbf{1}-\mathbf{G}(t)\right) )}}\geq {\frac{
\lambda(\delta^{2-\gamma-\beta}-1)}{(\delta-1)\delta(2-\gamma-\beta)}}.
\end{equation*}
Letting now $\delta\downarrow1$ yields
\begin{equation*}
\underset{t\to\infty}{\lim\inf}\,{\frac{(\mathbf{v,\Phi }(\mathbf{Q}(t;
\mathbf{1}-\lambda \mathbf{u}\psi (t))))}{\psi(t)(\mathbf{v}, \mathbf{u}
\otimes \left( \mathbf{1}-\mathbf{G}(t)\right) )}}\geq \lambda.
\end{equation*}
The proof of the converse inequality for the upper limit proceeds similarly,
starting with $\int_{\delta T}^T t(\mathbf{v,\Phi }(\mathbf{Q}(t;\mathbf{1}
-\lambda \mathbf{u}\psi (t)))) \mathrm{d}t$ for $\delta\in (0,1)$.

Using \eqref{inter9}, \eqref{DDEF} and Lemma \ref{L_ratio} we infer
\begin{equation*}
\lambda \psi (t)v_{2}u_{2}\left( 1-G_{2}(t)\right) \sim B{\mathfrak{Q}}
^{2}(t,\mathbf{1}-\lambda \mathbf{u}\psi (t))
\end{equation*}
and then
\begin{equation}
\mathbf{Q}(t;\mathbf{1}-\lambda \mathbf{u}\psi (t))\sim \sqrt{\lambda }\sqrt{%
\frac{v_{2}u_{2}}{B}\psi (t)(1-G_2(t)) } \mathbf{u}^{\mathbf{\dagger }}
\label{M2}
\end{equation}
by another appeal to Lemma \ref{L_ratio}. Since $\lim_{t\to\infty}\psi(t)=0$
and, in view of \eqref{Nzone1}, $\lim_{t\to\infty}\psi(t)(1-G_2(t))=0$ Lemma %
\ref{LRed} implies that \eqref{0M2} is a consequence of \eqref{M2}.

It remains to prove \eqref{Fpart} under additional assumption \eqref{M1}. To
this end, we first note that \eqref{0M2} and \eqref{Nzone1} together imply
\begin{equation*}
\lim_{N,t\rightarrow \infty }N\mathbf{Q}(t;\mathbf{e}^{-\lambda \mathbf{u}%
\psi (t)}) =r\sqrt{\lambda }\mathbf{u}^{\mathbf{\dagger }}.
\end{equation*}
Observe further that
\begin{eqnarray*}
\mathbf{0} &\leq &\mathbf{Q}(t;\mathbf{e}^{-\lambda \mathbf{u}\psi (t)})-%
\mathbf{Q}(t;1,e^{-\lambda u_{2}\psi (t)}) \\
&\leq &\lambda u_{1}\psi (t)\mathbf{P}_{1}\left( t\right) \sim \lambda \mu
_{1}\beta \Gamma _{\beta }u_{1}\frac{\psi (t)}{\mu _{2}(t)}\left(
D_{11},D_{21}\right) ^{\dag },
\end{eqnarray*}%
where the last equivalence follows from Lemma \ref{L_firsrMom}, which in
combination with \eqref{M1} yields
\begin{equation*}
N\mathbf{Q}(t;\mathbf{e}^{-\lambda \mathbf{u}\psi (t)})\sim N\mathbf{Q}%
(t;1,e^{-\lambda u_{2}\psi (t)}).
\end{equation*}
The proof of Lemma \ref{L_z2} is complete.

Passing to the proof of Theorem \ref{T_zone2} we note that (\ref{inter3})
follows from \eqref{Fpart} and (\ref{Blarge}), while \eqref{inter4} is a
consequence of \eqref{FinFirst1}, \eqref{Fpart}, \eqref{Blarge} and
\begin{eqnarray}
\mathbf{0} &\leq &\mathbf{Q}(t;0,e^{-\lambda u_{2}\psi (t)})-\mathbf{Q}%
(t;1,e^{-\lambda u_{2}\psi (t)})  \notag \\
&\leq &\mathbf{P}_{1}\left( t\right) \sim \mu _{1}\beta \Gamma _{\beta }%
\frac{1}{\mu _{2}(t)}\left( D_{11},D_{21}\right) ^{\dag}.  \label{inter6}
\end{eqnarray}
The proof of Theorem \ref{T_zone2} is complete.

\end{document}